\documentclass[11pt]{article}
\usepackage{txfonts}
\usepackage{color,epsfig}
\usepackage{mathrsfs}
\usepackage{amssymb}
\usepackage{amsfonts,graphicx,psfrag}

\newtheorem{Th}{Theorem}[section]

\newtheorem{lemma}{Lemma}[section]
\newtheorem{theorem}[lemma]{Theorem}
\newtheorem{corollary}[lemma]{Corollary}

\newtheorem{proposition}[lemma]{Proposition}
\newtheorem{definition}[lemma]{Definition}

\newtheorem{remark}[lemma]{Remark}
\let\lutzremark=\remark
\def\remark{\lutzremark\normalfont}

\newtheorem{notation}[lemma]{Notation}

\def\be{\begin{equation}}
\def\ee{\end{equation}}
\def\bea{\begin{eqnarray}}
\def\eea{\end{eqnarray}}
\def\bes{\begin{eqnarray*}}
\def\ees{\end{eqnarray*}}

\def\nn{\nonumber}
\def\<{\langle}
\def\>{\rangle}
\def\lb{\label}
\def\bs{\setminus}
\def\pt{\partial}

\def\R{{\bf R}}
\def\C{{\bf C}}
\def\Z{{\bf Z}}
\def\N{{\bf N}}
\def\U{{\bf U}}

\def\aa{{\alpha}}
\def\bb{{\beta}}
\def\ga{{\gamma}}
\def\Ga{{\Gamma}}

\def\th{{\theta}}
\def\Th{{\Theta}}
\def\om{{\omega}}
\def\Om{{\Omega}}
\def\ep{{\epsilon}}
\def\lm{{\lambda}}

\def\dl{{\delta}}

\def\sg{{\sigma}}

\def\P{{\cal P}}

\def\diag{{\rm diag}}
\def\constant{{\rm constant}}

\def\span{{\rm span}}
\def\Sp{{\rm Sp}}

\def\dm{{\rm \diamond}}

\def\hb{\vrule height0.18cm width0.14cm $\,$}
\def\ol#1{\overline{#1}}
\def\td#1{\tilde{#1}}

\title{Maslov-type indices and linear stability of elliptic Euler solutions of
the three-body problem}

\author{Qinglong Zhou$^{1} $\thanks{Partially supported by NSFC (No.11425105) of China.
           E-mail:zhouqinglong@sdu.edu.cn}\quad
        Yiming Long$^{2} $\thanks{Partially supported by NNSF, MCME, LPMC of MOE of China,
           Nankai University, and the Beijing Center for Math. and Info. Interdisciplinary Sciences
           at Capital Normal University. E-mail: longym@nankai.edu.cn}, \\
\\$^{1}$ School of Mathematics\\Shandong University, Jinan 250100, Shandong, China\\
  $^{2}$ Chern Institute of Mathematics and LPMC\\Nankai University, Tianjin 300071, China\\
}

\begin{document}

\maketitle

\begin{abstract}
{In this paper, we use the central configuration coordinate decomposition to study the linearized
Hamiltonian system near the $3$-body elliptic Euler solutions. Then using the Maslov-type $\omega$-index
theory of symplectic paths and the theory of linear operators we compute the $\om$-indices and obtain
certain properties of linear stability of the Euler elliptic solutions of the classical three-body
problem.}
\end{abstract}



\renewcommand{\theequation}{\thesection.\arabic{equation}}

\setcounter{equation}{0}\setcounter{figure}{0}
\section{Introduction and main results}
\label{sec:1}

In 1767, Euler (\cite{Euler}) discovered some celebrated periodic solutions, now named after him,
to the planar three-body problem, namely the three bodies are collinear at any instant
of the motion and at the same time each body travels along a specific Keplerian elliptic orbit about the
center of masses of the system. All these orbits are homographic solutions.
When $0\le e<1$, the Keplerian orbit is elliptic, we call such elliptic Euler (Lagrangian) solutions
Euler (Lagrangian) {\it elliptic relative equilibria}. Specially when $e=0$, the Keplerian elliptic
motion becomes circular motion and then all the three bodies move around the center of masses along circular
orbits with the same frequency, which are called Euler (Lagrangian) {\it relative equilibria} traditionally.
In this paper, we study the Maslov-type and Morse indices of such elliptic Euler solutions which are closely
related to their linear stability.

Denote by $q_1,q_2,q_3\in \R^2$ the position vectors of three particles with masses $m_1,m_2,m_3>0$
respectively. Then the system of equations for this problem is
\be   m_i\ddot{q}_i=\frac{\partial U}{\partial q_i}, \qquad {\rm for}\quad i=1, 2, 3, \lb{1.1}\ee
where $U(q)=U(q_1,q_2,q_3)=\sum_{1\leq i<j\leq 3}\frac{m_im_j}{|q_i-q_j|}$ is the
potential or force function by using the standard norm $|\cdot|$ of vector in $\R^2$.

Note that $2\pi$-periodic solutions of this problem correspond to critical points of the action functional
$$ \mathcal{A}(q)=\int_{0}^{2\pi}\left[\sum_{i=1}^3\frac{m_i|\dot{q}_i(t)|^2}{2}+U(q(t))\right]dt $$
defined on the loop space $W^{1,2}(\R/2\pi\Z,\hat{\mathcal {X}})$, where
$$  \hat{\mathcal {X}}:=\left\{q=(q_1,q_2,q_3)\in (\R^2)^3\,\,\left|\,\,
       \sum_{i=1}^3 m_iq_i=0,\,\,q_i\neq q_j,\,\,\forall i\neq j \right. \right\}  $$
is the configuration space of the planar three-body problem.

Letting $p_i=m_i\dot{q}_i\in\R^2$ for $1\le i\le 3$, then (\ref{1.1}) is transformed to a Hamiltonian system
\be \dot{p}_i=-\frac{\partial H}{\partial q_i},\,\,\dot{q}_i
  = \frac{\partial H}{\partial p_i},\qquad {\rm for}\quad i=1,2,3,  \lb{1.2}\ee
with Hamiltonian function
\be H(p,q)=H(p_1,p_2,p_3, q_1,q_2,q_3)=\sum_{i=1}^3\frac{|p_i|^2}{2m_i}-U(q_1,q_2,q_3).  \lb{1.3}\ee

For the planar three-body problem with masses $m_1, m_2, m_3>0$, it turns out that the stability
of elliptic Euler solutions depends on two parameters, namely the mass parameter $\beta\in [0,7]$
defined below and the eccentricity $e\in [0,1)$,
\begin{equation}
\beta=\frac{m_1(3x^2+3x+1)+m_3x^2(x^2+3x+3)}{x^2+m_2[(x+1)^2(x^2+1)-x^2]},    \lb{1.4}
\end{equation}
where $x$ is the unique positive solution of the Euler quintic polynomial equation (\ref{quintic.polynomial}).

The linear stability of Lagrangian relative equilibria can be found in
Gascheau (\cite{Ga}, 1843), Routh (\cite{R2}, 1875), Danby (\cite{Dan}, 1964) and Roberts (\cite{R1}, 2002).
In 2005, Meyer and Schmidt (cf. \cite{MS}) used heavily the central configuration nature
of the elliptic Lagrangian orbits and decomposed the fundamental solution of the elliptic Lagrangian
orbit into two parts symplectically, one of which is the same as that of the Keplerian solution and
the other is the essential part for the stability.

In 2004-2006, Mart\'{\i}nez, Sam\`{a} and Sim\'{o} (\cite{MSS},\cite{MSS1},\cite{MSS2})
studied the stability problem including Euler elliptic relative equilibria
when $e>0$ is small enough by using normal form theory, and $e<1$
and close to $1$ enough by using blow-up technique in general homogeneous potential. They
further gave a much more complete bifurcation diagram numerically and a beautiful figure was
drawn there for the full $(\bb,e)$ range (cf. Figure 4 of \cite{MSS2}).

In \cite{HS1} and \cite{HS2} of 2009-2010, Hu and Sun found a new way to relate the stability problem to the
iterated Morse indices. Recently, by observing new phenomenons and discovering new properties of elliptic
Lagrangian solution, in the joint paper \cite{HLS} of Hu, Long and Sun, the linear stability of elliptic
Lagrangian solution is completely solved analytically by index theory (cf. \cite{Lon1} and \cite{Lon4}) and
the new results are related directly to $(\beta,e)$ in the full parameter rectangle.

In the current paper, for the elliptic Euler solutions, following the central configuration coordinate method of
Meyer and Schmidt in \cite{MS} and the index method used by Hu, Long and Sun in \cite{HLS}, we linearized the
Hamiltonian system (\ref{1.2})-(\ref{1.3}) near the Euler elliptic solution in Section 2 below. Here the linearized
Hamiltonian system can also be decomposed into two parts symplectically, one of which is the same as that of the
Kepler solutions, and the other is a $4$-dimensional Hamiltonian system whose fundamental solution is the essential
part for the stability of the elliptic Euler solutions. However, the essential part here is very different from
that of the Lagrangian elliptic solutions in \cite{MS} and \cite{HLS}. This essential part is denoted by
$\gamma_{\beta,e}(t)$ for $t\in [0,2\pi]$, which is a path in $\Sp(4)$ starting from the identity. Then we use
index theory to compute the Maslov-type indices of $\ga_{\bb,e}$ and determine its stability properties.

Following \cite{Lon2} and \cite{Lon4}, for any $\omega\in\U=\{z\in\C\;|\;|z|=1\}$ we can define a real function
$D_\om(M)=(-1)^{n-1}\overline{\om}^n det(M-\om I_{2n})$ for any $M$ in the symplectic group $\Sp(2n)$.
Then we can define $\Sp(2n)_{\om}^0 = \{M\in\Sp(2n)\,|\, D_{\om}(M)=0\}$ and
$\Sp(2n)_{\om}^{\ast} = \Sp(2n)\bs \Sp(2n)_{\om}^0$. The orientation of $\Sp(2n)_{\om}^0$ at any of its point
$M$ is defined to be the positive direction $\frac{d}{dt}Me^{t J}|_{t=0}$ of the path $Me^{t J}$ with $t>0$ small
enough. Let $\nu_{\om}(M)=\dim_{\C}\ker_{\C}(M-\om I_{2n})$. Let
$\mathcal{P}_{2\pi}(2n) = \{\ga\in C([0,2\pi],\Sp(2n))\;|\;\ga(0)=I\}$ and
$\xi(t)=\diag(2-\frac{t}{2\pi}, (2-\frac{t}{2\pi})^{-1})$ for $0\le t\le 2\pi$.

Given any two $2m_k\times 2m_k$ matrices of square block form
$M_k=\left(\matrix{A_k&B_k\cr
                   C_k&D_k\cr}\right)$ with $k=1, 2$,
the symplectic sum of $M_1$ and $M_2$ is defined (cf. \cite{Lon2} and \cite{Lon4}) by
the following $2(m_1+m_2)\times 2(m_1+m_2)$ matrix $M_1\dm M_2$:
$$
M_1\dm M_2=\left(\matrix{A_1 &   0 & B_1 &   0\cr
                             0   & A_2 &   0 & B_2\cr
                             C_1 &   0 & D_1 &   0\cr
                             0   & C_2 &   0 & D_2\cr}\right),
$$
and $M^{\dm k}$ denotes the $k$ copy $\dm$-sum of $M$. For any two paths $\ga_j\in\P_{\tau}(2n_j)$
with $j=0$ and $1$, let $\ga_0\dm\ga_1(t)= \ga_0(t)\dm\ga_1(t)$ for all $t\in [0,\tau]$.

For any $\ga\in \mathcal{P}_{2\pi}(2n)$ we define $\nu_\om(\ga)=\nu_\om(\ga(2\pi))$ and
$$  i_\om(\ga)=[\Sp(2n)_\om^0:\ga\ast\xi^n], \qquad {\rm if}\;\;\ga(2\pi)\not\in \Sp(2n)_{\om}^0,  $$
i.e., the usual homotopy intersection number, and the orientation of the joint path $\ga\ast\xi_n$ is
its positive time direction under homotopy with fixed end points. When $\ga(2\pi)\in \Sp(2n)_{\om}^0$,
we define $i_{\om}(\ga)$ be the index of the left rotation perturbation path $\ga_{-\ep}$ with $\ep>0$
small enough (cf. Def. 5.4.2 on p.129 of \cite{Lon4}). The pair
$(i_{\om}(\ga), \nu_{\om}(\ga)) \in \Z\times \{0,1,\ldots,2n\}$ is called the index function of $\ga$
at $\om$. When $\nu_{\om}(\ga)=0$ or $\nu_{\om}(\ga)>0$, the path $\ga$ is called
$\om$-{\it non-degenerate} or $\om$-{\it degenerate} respectively. For more details we refer to the
Appendix 5.2 or \cite{Lon4}.

The following three theorems describe main results proved in this paper.

\begin{theorem}\label{T1.1}
In the planar three-body problem with masses $m_1, m_2$, and $m_3>0$, for the elliptic Euler solution
$q=q_{\bb,e}(t)=(q_1(t), q_2(t), q_3(t))$ with eccentricity $e$ and mass parameter $\beta$ given by
(\ref{1.4}), we denote by $\ga_{\bb,e}:[0,2\pi] \to \Sp(4)$ the essential part of the fundamental solution
of the linearized Hamiltonian system of (\ref{1.1}) at $q$. Then the following results on the Maslov-type
indices of $\ga_{\bb,e}$ hold.

(i) $(i_1(\ga_{0,e}), \nu_1(\ga_{0,e})) = (0,3)$ and $(i_{\om}(\ga_{0,e}), \nu_{\om}(\ga_{0,e})) = (2,0)$
for $\om\in\U\bs\{1\}$.

(ii) Let
\be  \hat\bb_n=\frac{n^2-3+\sqrt{9n^4-14n^2+9}}{4} \qquad \forall\;n\in\N.  \lb{1.5}\ee
Then
\begin{eqnarray}
&& i_1(\ga_{\bb,0}) = \left\{\matrix{
   0, &  {\it if}\;\;\bb=\hat\bb_1=0, \cr
   2n+1, &  {\it if}\;\;\bb\in(\hat\bb_n,\hat\bb_{n+1}]\;\;{\it for}\;n\in\N, \cr}\right. \lb{1.6}\\
&& \nu_1(\ga_{\bb,0}) = \left\{\matrix{
   3, &  {\it if}\;\;\bb=\hat\bb_1=0, \cr
   2, &  {\it if}\;\;\bb=\hat\bb_n,\;n\ge2, \cr
   0, &  {\it if}\;\;\bb\in [0,+\infty)\bs\{\hat\bb_n\;|\;n\in\N\}. \cr}\right. \lb{1.7}
\end{eqnarray}

(iii) Let
\be   \hat\bb_{n+\frac{1}{2}}=\frac{(n+\frac{1}{2})^2-3+\sqrt{9(n+\frac{1}{2})^4-14(n+\frac{1}{2})^2+9}}{4}
        \qquad \forall\;n\in\N.  \lb{1.8}\ee
Then
\begin{eqnarray}
&& i_{-1}(\ga_{\bb,0}) = \left\{\matrix{
   2, &  {\it if}\;\;\bb\in[0,\hat\bb_{\frac{3}{2}}], \cr
   2n, &  {\it if}\;\;\bb\in(\hat\bb_{n-\frac{1}{2}},\hat\bb_{n+\frac{1}{2}}]\;\;{\it for}\;n\ge2, \cr}\right. \lb{1.9}\\
&& \nu_{-1}(\ga_{\bb,0}) = \left\{\matrix{
   2, &  {\rm if}\;\;\bb=\hat\bb_{n+\frac{1}{2}}\;{\it for}\;n\in\N, \cr
   0, &  {\rm if}\;\;\bb\in [0,+\infty)\bs\{\hat\bb_{n+\frac{1}{2}}\;|\;n\in\N\}.\cr} \right.
                                   \lb{1.10}
\end{eqnarray}

(iv) For fixed $e\in [0,1)$ and $\om\in\U$, $i_{\om}(\ga_{\bb,e})$ is non-decreasing and tends to $+\infty$
when $\bb$ increases from $0$ to $+\infty$.

(v) $i_1(\ga_{\bb,e})>0$ is odd for all $(\bb,e)\in (0,+\infty)\times [0,1)$.

(vi) $i_1(\ga_{\bb,e})\le 4n+2$ holds when $\bb<\frac{2}{3\sqrt{2}-1}(n^2-\frac{e}{1+e})(1-e)-1$ for
any $n\in\N$.

(vii) For any $\epsilon\in(0,1)$, there exists a $\bb_\epsilon>0$ such that $\nu_1(\ga_{\bb,e})=0$, i.e.,
$\ga_{\bb,e}$ is non-degenerate when $(\bb,e)\in(0,\bb_\epsilon]\times[0,1-\epsilon]$.
\end{theorem}

\begin{remark}\label{R1.2}
(i) Here we are specially interested in indices in eigenvalues $1$ and $-1$. The reason is that the
major changes of the linear stability of the elliptic Euler solutions happen near the eigenvalues
$1$ and $-1$, and such information is used in the next theorem to get the separation curves of the
linear stability domain $[0,+\infty)\times [0,1)$ of the mass and eccentricity parameter $(\bb,e)$.

(ii) The situations of other eigenvalues $\om\in \U\bs\R$ of $\ga_{\bb,e}(2\pi)$ can be obtained by
the method in Section 4 below similarly, which then yields complete understanding on the eigenvalue
distribution of $\ga_{\bb,0}(2\pi)$ for all $\bb\ge 0$, i.e., the linear stability of the Euler relative
equilibria $q_{\bb,0}(t)$. Note that by the essential part of the linearized Hamiltonian system at
the elliptic Euler solutions found in (\ref{2.19}) below, $e=0$ yields an autonomous Hamiltonian
system, and thus the linear stability is explicitly computable.

(iii) Note that $\bb\in [0,7]$ in its physical meaning. For mathematical interest and convenience,
we extend the range of the parameter $\bb$ to $[0,\infty)$.
\end{remark}

\begin{theorem}\label{T1.3}
Using notations in Theorem \ref{T1.1}, for the elliptic Euler solution $q=q_{\bb,e}(t)$ with
eccentricity $e$ and mass parameter $\beta$ given by (\ref{1.4}), the following results on the
linear stability separation curves of $\ga_{\bb,e}$ in the parameter $(\bb,e)$ domain
$\Th=[0,+\infty)\times [0,1)$ hold.
Letting
\bea
\Ga_n &=& \{(\bb_{2n-1}(1,e),e)\;|\;e\in [0,1)\}\quad {\it with}\quad \bb_{2n-1}(1,e)=(\bb_{2n}(1,e), \nn\\
\Xi_n^- &=& \{(\bb_{2n-1}(-1,e),e)\;|\;e\in [0,1)\},  \nn\\
\Xi_n^+ &=& \{(\bb_{2n}(-1,e),e)\;|\;e\in [0,1)\}, \nn\eea
we then have the following:

(i) Starting from the point $(\hat\bb_{n+1},0)$ defined in (\ref{1.5}) for $n\in\N$, there exists
exactly one $1$-degenerate curve $\Ga_n$ of $\ga_{\bb,e}(2\pi)$ which is perpendicular to the
$\bb$-axis, goes up into the domain $\Th$, intersects each horizontal line $e=\constant$ in $\Th$
precisely once for each $e\in (0,1)$, and satisfies $\nu_1(\ga_{\bb_{2n}(1,e),e})=2$ at such an
intersection point $(\bb_{2n}(1,e),e)\in\Ga_n$, see Figure 1 below (cf. left figure of Figure 6 in
\cite{MSS1}). Further more, $\bb_{2n}(1,e)$ is a real analytic function in $e\in [0,1)$.

(ii) Starting from the point $(\hat\bb_{n+1/2},0)$ defined in (\ref{1.8}) for $n\in\N$, there exists
exactly two $-1$-degenerate curves $\Xi_n^{\pm}$ of $\ga_{\bb,e}(2\pi)$ which are perpendicular
to the $\bb$-axis, go up into the domain $\Th$. Moreover, for each $e\in (0,1)$,
if $\bb_{2n-1}(-1,e)\ne\bb_{2n}(-1,e)$, the two curves intersect each horizontal line $e=\constant$ in $\Th$
precisely once and satisfy $\nu_1(\ga_{\bb_{2n-1}(-1,e),e})=\nu_1(\ga_{\bb_{2n}(-1,e),e})=1$ at such an
intersection point $(\bb_{2n-1}(-1,e),e)\in\Xi_n^{-}$ and $(\bb_{2n}(-1,e),e)\in\Xi_n^{+}$;
if $\bb_{2n-1}(-1,e)=\bb_{2n}(-1,e)$, the two curves intersect each horizontal line $e=\constant$ in $\Th$
at the same point and satisfy $\nu_1(\ga_{\bb_{2n-1}(-1,e),e})=2$ at such an intersection point
$(\bb_{2n-1}(-1,e),e)\in\Xi_n^{+}\cap\Xi_n^{-}$. Further more, both $\bb_{2n-1}(-1,e)$ and
$\bb_{2n}(-1,e)$ are real piecewise analytic functions in $e\in [0,1)$. Note that in Figure 1 below the two curves
which start from the point $(\hat\bb_{n+1/2},0)$ where $n\ge2$ are close enough, so they look like just
one curve in our figure.

(iii) The $1$-degenerate curves and $-1$-degenerate curves of the
elliptic Euler solutions in Figure 1 can be ordered from left to right by
\begin{equation}
0,\; \Xi_1^-,\; \Xi_1^+,\; \Ga_1,\; \Xi_2^-,\; \Xi_2^+,\; \Ga_2,\; \ldots,\; \Xi_n^-,\; \Xi_n^+,\; \Ga_n,\;\ldots .
\end{equation}
Moreover, for $n_1,n_2\in\N$, $\Ga_{n_1}$ and $\Xi^{\pm}_{n_2}$ cannot intersect each other;
if $n_1\ne n_2$, $\Ga_{n_1}$ and $\Ga_{n_2}$ cannot intersect each other,
and $\Xi^{\pm}_{n_1}$ and $\Xi^{\pm}_{n_2}$ cannot intersect each other.
More precisely, for each fixed $e\in [0,1)$, we have
\begin{eqnarray}
0&<&\bb_1(-1,e)\le \bb_2(-1,e)<\bb_1(1,e)=\bb_2(1,e)<\bb_3(-1,e)\le \bb_4(-1,e)<\bb_3(1,e)=\bb_4(1,e)<\cdots
\nonumber
\\
&<&\bb_{2n-1}(-1,e)\le \bb_{2n}(-1,e)<\bb_{2n-1}(1,e)=\bb_{2n}(1,e)<\cdots
\end{eqnarray}
\end{theorem}

\begin{remark}\label{R1.4}
We refer readers to the recent interesting paper \cite{HO2} of Professor Xijun Hu and Dr. Yuwei Ou, which
appeared almost simultaneously with the first version of the current paper \cite{ZhoL1}. In \cite{HO2} the
authors introduced the collision index, studied the behavior of the above $1$-degenerate and $-1$-degenerate
curves as $e\to 1$, and completely understood the properties of these curves when $e$ is close to $1$.
Note that our Theorems 1.1, 1.3 and 1.5 below together with the results in \cite{HO2} give a complete
analytical understanding of the stability properties of the $3$-body elliptic Euler solutions.
\end{remark}

\begin{figure}[ht]
\centering
\vskip 1 cm
\includegraphics[height=8.5cm]{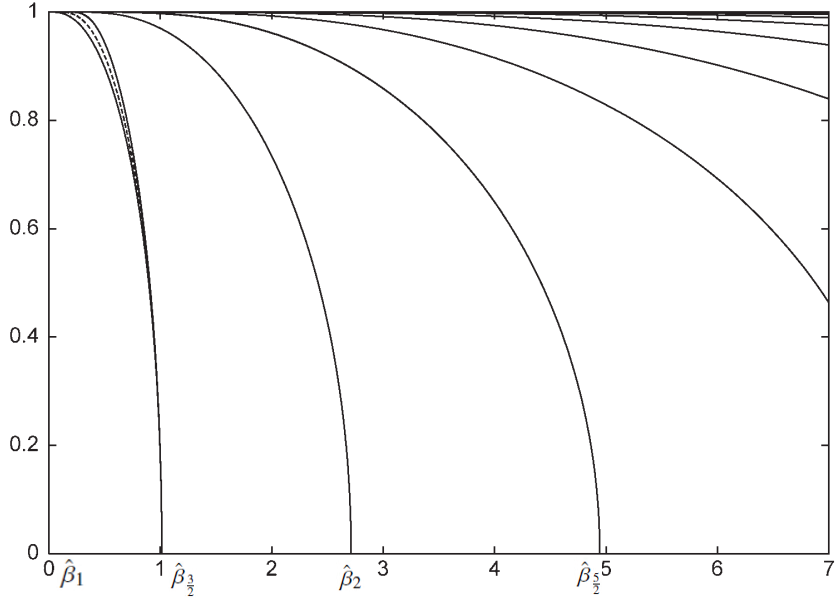}
\caption{The $1$-degenerate and $-1$-degenerate curves of Euler elliptic relative equilibria of
the planar three-body problem in the $(\beta,e)$ rectangle $[0,7]\times [0,1)$.}
\end{figure}

The concept of ``$M\approx N$" for two symplectic matrices $M$ and $N$, i.e., $N\in \Omega^0(M)$,
was first introduced in \cite{Lon2} of 1999, which can be found in the Definition 5.2 of the
Appendix 5.2 in this paper following Definition 1.8.5 of \cite{Lon4}. This notion is broader than
the symplectic similarity in general as pointed out on p.38 of \cite{Lon4}.

For the normal forms of $\ga_{\bb,e}(2\pi)$, we have the following theorem.

\begin{theorem}\label{T1.5}
For the normal forms of $\ga_{\bb,e}(2\pi)$ when $\bb\ge0,0\le e<1$, for $n\in\N$,
we have the following results:

(i) If $\bb=0$, we have $i_1(\ga_{0,e}(2\pi))=0$, $\nu_1(\ga_{0,e}(2\pi))=3$,
$i_{-1}(\ga_{0,e}(2\pi))=2$, $\nu_{-1}(\ga_{0,e}(2\pi))=0$
and $\ga_{0,e}(2\pi)\approx I_2\diamond N_1(1,1)$;

(ii) If $0<\bb<\bb_1(-1,e)$, we have $i_1(\ga_{\bb,e}(2\pi))=3$, $\nu_1(\ga_{\bb,e}(2\pi))=0$,
$i_{-1}(\ga_{\bb,e}(2\pi))=2$, $\nu_{-1}(\ga_{\bb,e}(2\pi))=0$
and $\ga_{\bb,e}(2\pi)\approx R(\theta)\diamond D(2)$ for some $\theta\in(0,\pi)$;

(iii) If $\bb=\bb_1(-1,e)=\bb_2(-1,e)$, we have $i_1(\ga_{\bb,e}(2\pi))=3$, $\nu_1(\ga_{\bb,e}(2\pi))=0$,
$i_{-1}(\ga_{\bb,e}(2\pi))=2$, $\nu_{-1}(\ga_{\bb,e}(2\pi))=2$
and $\ga_{\bb,e}(2\pi)\approx -I_2\diamond D(2)$;

(iv) If $\bb_1(-1,e)\ne \bb_2(-1,e)$ and $\bb=\bb_1(-1,e)$,
we have $i_1(\ga_{\bb,e}(2\pi))=3$, $\nu_1(\ga_{\bb,e}(2\pi))=0$,
$i_{-1}(\ga_{\bb,e}(2\pi))=2$, $\nu_{-1}(\ga_{\bb,e}(2\pi))=1$
and $\ga_{\bb,e}(2\pi)\approx N_1(-1,-1)\diamond D(2)$;

(v) If $\bb_1(-1,e)\ne \bb_2(-1,e)$ and $\bb_1(-1,e)<\bb<\bb_2(-1,e)$,
we have $i_1(\ga_{\bb,e}(2\pi))=3$, $\nu_1(\ga_{\bb,e}(2\pi))=0$,
$i_{-1}(\ga_{\bb,e}(2\pi))=3$, $\nu_{-1}(\ga_{\bb,e}(2\pi))=0$
and $\ga_{\bb,e}(2\pi)\approx D(-2)\diamond D(2)$;

(vi) If $\bb_1(-1,e)\ne \bb_2(-1,e)$ and $\bb=\bb_2(-1,e)$,
we have $i_1(\ga_{\bb,e}(2\pi))=3$, $\nu_1(\ga_{\bb,e}(2\pi))=0$,
$i_{-1}(\ga_{\bb,e}(2\pi))=3$, $\nu_{-1}(\ga_{\bb,e}(2\pi))=1$
and $\ga_{\bb,e}(2\pi)\approx N_1(-1,1)\diamond D(2)$;

(vii) If $\bb_2(-1,e)<\bb<\bb_1(1,e)$,
we have $i_1(\ga_{\bb,e}(2\pi))=3$, $\nu_1(\ga_{\bb,e}(2\pi))=0$,
$i_{-1}(\ga_{\bb,e}(2\pi))=4$, $\nu_{-1}(\ga_{\bb,e}(2\pi))=0$
and $\ga_{\bb,e}(2\pi)\approx R(\theta)\diamond D(2)$
for some $\theta\in(\pi,2\pi)$;

(viii) If $\bb=\bb_{2n-1}(1,e)(=\bb_{2n}(1,e))$,
we have $i_1(\ga_{\bb,e}(2\pi))=2n+1$, $\nu_1(\ga_{\bb,e}(2\pi))=2$,
$i_{-1}(\ga_{\bb,e}(2\pi))=2n+2$, $\nu_{-1}(\ga_{\bb,e}(2\pi))=0$
and $\ga_{\bb,e}(2\pi)\approx I_2\diamond D(2)$;

(ix) If $\bb_{2n}(1,e)<\bb<\bb_{2n+1}(-1,e)$,
we have $i_1(\ga_{\bb,e}(2\pi))=2n+3$, $\nu_1(\ga_{\bb,e}(2\pi))=0$,
$i_{-1}(\ga_{\bb,e}(2\pi))=2n+2$, $\nu_{-1}(\ga_{\bb,e}(2\pi))=0$
and $\ga_{\bb,e}(2\pi)\approx R(\theta)\diamond D(2)$ for some $\theta\in(0,\pi)$;

(x) If $\bb=\bb_{2n+1}(-1,e)=\bb_{2n+2}(-1,e)$,
we have $i_1(\ga_{\bb,e}(2\pi))=2n+3$, $\nu_1(\ga_{\bb,e}(2\pi))=0$,
$i_{-1}(\ga_{\bb,e}(2\pi))=2n+2$, $\nu_{-1}(\ga_{\bb,e}(2\pi))=2$
and $\ga_{\bb,e}(2\pi)\approx -I_2\diamond D(2)$;

(xi) If $\bb_{2n+1}(-1,e)\ne \bb_{2n+2}(-1,e)$ and $\bb=\bb_{2n+1}(-1,e)$,
we have $i_1(\ga_{\bb,e}(2\pi))=2n+3$, $\nu_1(\ga_{\bb,e}(2\pi))=0$,
$i_{-1}(\ga_{\bb,e}(2\pi))=2n+2$, $\nu_{-1}(\ga_{\bb,e}(2\pi))=1$
and $\ga_{\bb,e}(2\pi)\approx N_1(-1,-1)\diamond D(2)$;

(xii) If $\bb_{2n+1}(-1,e)\ne \bb_{2n+2}(-1,e)$ and $\bb_{2n+1}(-1,e)<\bb<\bb_{2n+2}(-1,e)$,
we have $i_1(\ga_{\bb,e}(2\pi))=2n+3$, $\nu_1(\ga_{\bb,e}(2\pi))=0$,
$i_{-1}(\ga_{\bb,e}(2\pi))=2n+3$, $\nu_{-1}(\ga_{\bb,e}(2\pi))=0$
and $\ga_{\bb,e}(2\pi)\approx D(-2)\diamond D(2)$;

(xiii) If $\bb_{2n+1}(-1,e)\ne \bb_{2n+2}(-1,e)$ and $\bb=\bb_{2n+2}(-1,e)$,
we have $i_1(\ga_{\bb,e}(2\pi))=2n+3$, $\nu_1(\ga_{\bb,e}(2\pi))=0$,
$i_{-1}(\ga_{\bb,e}(2\pi))=2n+3$, $\nu_{-1}(\ga_{\bb,e}(2\pi))=1$
and $\ga_{\bb,e}(2\pi)\approx N_1(-1,1)\diamond D(2)$;

(xiv) If $\bb_{2n+2}(-1,e)<\bb<\bb_{2n+1}(1,e)$,
we have $i_1(\ga_{\bb,e}(2\pi))=2n+3$, $\nu_1(\ga_{\bb,e}(2\pi))=0$,
$i_{-1}(\ga_{\bb,e}(2\pi))=2n+4$, $\nu_{-1}(\ga_{\bb,e}(2\pi))=0$
and $\ga_{\bb,e}(2\pi)\approx R(\theta)\diamond D(2)$
for some $\theta\in(\pi,2\pi)$.
\end{theorem}

In the proof of these theorems, motivated by the techniques of \cite{HLS}, we study properties of the symplectic path
$\ga_{\bb,e}$ in $\Sp(4)$ and the second order differential operators $A(\bb,e)$ corresponding to $\ga_{\bb,e}$.
To get the information on the indices of $\ga_{\bb,e}$ for $(\bb,e)\in \Th$, one of the main ingredients of the
proof is the non-decreasing property of $\om$-index proved in Lemma \ref{Lemma:increasing.of.index} and Corollary \ref{C4.5} below
for all $\om\in\U$.

The rest of this paper is focused on the proof of Theorems 1.1, 1.3 and 1.5. For Theorem 1.1, the index properties in
(i)-(iii) are established in Section 3;
the non-decreasing property (iv) is proved in Corollary \ref{C4.5};
the property (v)
is proved in Theorem \ref{Th:odd.indices}; the estimate (vi) is proved in Proposition 4.4; and the non-degenerate property (vii)
is proved in Theorem 4.6.
Theorem 1.5 is proved in the Subsection 4.3.
For Theorem 1.3, (i) on the $1$-degenerate curves $\Ga_n$ is proved in Subsection 4.3 and Subsection 4.4;
(ii) on the $-1$-degenerate curves $\Xi_n$ is proved in the Subsection 4.4; and (iii) is prove in the Subsection 4.3.

\setcounter{equation}{0}
\section{Preliminaries}\label{sec:2} 

In the subsection 5.2 of the Appendix, we give a brief review on the Maslov-type $\om$-index theory for
$\om$ in the unit circle of the complex plane following \cite{Lon4}. In the following, we use notations
introduced there.

\subsection{The essential part of the fundamental solution of the elliptic Euler orbit}

In \cite{MS} (cf. p.275), Meyer and Schmidt gave the essential part of the fundamental solution of the
elliptic Lagrangian orbit. Their method is explained in \cite{Lon5} too. Our study on elliptic Euler
solutions is based upon their method.

Suppose the three particles are all on the $x$-axis, $q_1=0$, $q_2=(x\alpha,0)^T$ and $q_3=((1+x)\alpha,0)^T$
for $\alpha=|q_2-q_3|>0$, $x\alpha=|q_1-q_2|$ and some $x>0$.
When $q_1$, $q_2$ and $q_3$ form a collinear central configurations, $x$ must satisfy Euler's quintic
equation as in p.148 of \cite{Euler}, p.276 of \cite{Win1} and p.29 of \cite{Lon5}:
\begin{equation}\label{quintic.polynomial}
(m_3+m_2)x^5+(3m_3+2m_2)x^4+(3m_3+m_2)x^3-(3m_1+m_2)x^2-(3m_1+2m_2)x-(m_1+m_2)=0.
\end{equation}
Moreover, by Descartes' rule of signs for polynomials (cf. p.300 of \cite{Jac1}),
polynomial (\ref{quintic.polynomial}) has only one positive solution $x$.

Without lose of generality, we normalize the three masses by
\begin{equation}\label{nomorlize.the.masses}
m_1+m_2+m_3=1.
\end{equation}
Then the center of mass of the three particles is
$$ q_0=m_1q_1+m_2q_2+m_3q_3=([m_2x+m_3(1+x)]\alpha,0)^T=([m_3+(1-m_1)x]\alpha,0)^T, $$
where we used (\ref{nomorlize.the.masses}) in the last equality.

For $i=1,2,3$, let $a_i=q_i-q_0$, and denote by $a_{ix}$ and $a_{iy}$ the $x$ and
$y$-coordinates of $a_i$ respectively. Then we have
\begin{equation}\label{a_ix}
a_{1x}=-[m_3+(1-m_1)x]\alpha,\quad
a_{2x}=(-m_3+m_1x)\alpha,\quad
a_{3x}=[(1-m_3)+m_1x]\alpha
\end{equation}
and
\begin{equation}\label{a_iy}
a_{iy}=0,\quad\quad {\rm for}\; i=1,2,3.
\end{equation}
Scaling $\alpha$ by setting $\sum_{i=1}^3m_i|a_i|^2=1$, we obtain
\begin{eqnarray}\label{alpha.square}
\alpha^2&=&\frac{\alpha^2}{\sum_{i=1}^3m_i|a_i|^2}=\frac{1}{m_1[-m_3-(1-m_1)x]^2+m_2[-m_3+m_1x]^2+m_3[1-m_3+m_1x]^2}
\nonumber\\
&=&\frac{1}{m_1(1-m_1)x^2+2m_1m_3x+m_3(1-m_3)}.
\end{eqnarray}

Now as in p.263 of \cite{MS}, Section 11.2 of \cite{Lon5}, we define
\begin{equation}\label{PQYX}
P=\left(\matrix{p_1\cr p_2\cr p_3}\right),
\quad
Q=\left(\matrix{q_1\cr q_2\cr q_3}\right),
\quad
Y=\left(\matrix{G\cr Z\cr W}\right),
\quad
X=\left(\matrix{g\cr z\cr w}\right),
\end{equation}
where $p_i$, $q_i$, $i=1,2,3$ and $G$, $Z$, $W$, $g$, $z$, $w$ are all column vectors in $\R^2$.
We make the symplectic coordinate change
\be\lb{transform1}  P=A^{-T}Y,\quad Q=AX,  \ee
where the matrix $A$ is constructed as in the proof of Proposition 2.1 in \cite{MS}.
Concretely, the matrix $A\in {\bf GL}(\R^6)$ is given by
\begin{equation}
A=
\left(
\matrix{
I\quad A_1\quad B_1\cr
I\quad A_2\quad B_2\cr
I\quad A_3\quad B_3
}
\right),
\end{equation}
where by (\ref{a_ix})-(\ref{a_iy}), each $A_i$ is a $2\times2$ matrix given by
\begin{equation}\label{Aa}
A_i = (a_i, Ja_i)=\left(\matrix{a_{ix}\quad 0\cr 0\quad a_{ix}}\right)=a_{ix}I.
\end{equation}
with $J=\left(\matrix{0& -1\cr 1& 0}\right)$.

To fulfill $A^TMA=I$ (cf. (13) in p.263 of \cite{MS}), we must have
\begin{eqnarray}
B_1&=&\rho_1(A_3-A_2)^T=\rho_1(a_{3x}-a_{2x})I=\rho_1\alpha I, \nn\\
B_2&=&\rho_2(A_1-A_3)^T=\rho_2(a_{1x}-a_{3x})I=-\rho_2(1+x)\alpha I, \nn\\
B_3&=&\rho_3(A_2-A_1)^T=\rho_3(a_{2x}-a_{1x})I=\rho_3x\alpha I, \nn
\end{eqnarray}
where
\begin{equation}\label{rho_i}
\rho_i=\frac{\sqrt{m_1m_2m_3}}{m_i},\quad \forall 1\le i\le 3.
\end{equation}
Denote by
\begin{equation}\label{b}
b_1=\rho_1\alpha,\quad b_2=-\rho_2(1+x)\alpha,\quad b_3=\rho_3x\alpha.
\end{equation}
Then we simply have
\begin{equation}\label{Bb}
B_i=b_iI,\quad \forall 1\le i\le 3.
\end{equation}

Under the coordinate change (\ref{transform1}), we get the kinetic energy
\begin{equation}
K=\frac{1}{2}(|G|^2+|Z|^2+|W|^2),
\end{equation}
and the potential function
\begin{equation}\label{U_ij}
U(z,w)=\sum_{1\le i<j\le3}U_{ij}(z,w),\quad U_{ij}(z,w)=\frac{m_im_j}{d_{ij}(z,w)},
\end{equation}
with
\begin{equation}
d_{ij}(z,w)=|(A_i-A_j)z+(B_i-B_j)w|=|(a_{ix}-a_{jx})z+(b_i-b_j)w|,
\end{equation}
where we used (\ref{Aa}) and (\ref{Bb}).

Let $\theta$ be the true anomaly. In \cite{MS}, Meyer and Schmidt introduced their celebrated
central configuration coordinates, which greatly simplified the corresponding systems. Then
under the same steps of symplectic transformation in the proof of Lemma 3.1 in \cite{MS}, the
resulting Hamiltonian function of the 3-body problem is given by
\be  H(\theta,\bar{Z},\bar{W},\bar{z},\bar{w})
= \frac{1}{2}(|\bar{Z}|^2+|\bar{W}|^2)
  + (\bar{z}\cdot J\bar{Z}+\bar{w}\cdot J\bar{W})+\frac{p-r(\theta)}{2p}(|\bar{z}|^2+|\bar{w}|^2)
  - \frac{r(\theta)}{\sigma}U(\bar{z},\bar{w}), \lb{new.H.function}\ee
where
\be  r(\theta)=\frac{p}{1+e\cos\theta},  \lb{r1}\ee
and
\be  \mu
= \sum_{1\le i<j\le3}\frac{m_im_j}{|a_i-a_j|}
= \frac{1}{\alpha}\left(\frac{m_1m_2}{x}+m_2m_3+\frac{m_3m_1}{1+x}\right), \quad
\sigma = (\mu p)^{1/4}.  \label{mu}\ee

Note that here as pointed out in Section 11 of \cite{Lon5}, the original constant $\sg=\mu p$ in
the line 9 on p.273 of \cite{MS} is not correct and should be corrected to $\sg=(\mu p)^{1/4}$.
Because this constant and the related corrections in this derivation are crucial in the later
computations of the linear stability, we refer readers to Section 2 of \cite{ZhoL2} for the
complete details of derivations of (\ref{new.H.function})-(\ref{mu}).

Indeed, $H$ given by (\ref{new.H.function}) is essentially the Hamiltonian of the system in the
pulsating frame, in which $\theta$ is the new independent variable, and $p=a(1-e^2)$ with
$a$ and $e$ being the semi-major axis and the eccentricity of $z(t)$ respectively.

We now derived the linearized Hamiltonian system at the Euler elliptic solutions.

\begin{proposition}\label{linearized.Hamiltonian}
Using notations in (\ref{PQYX}), elliptic Euler solution $(P(t),Q(t))^T$ of the system (\ref{1.2}) with
\begin{equation}
Q(t)=(r(t)R(\theta(t))a_1,r(t)R(\theta(t))a_2,r(t)R(\theta(t))a_3)^T,\quad P(t)=M\dot{Q}(t)
\end{equation}
in time $t$ with the matrix $M=diag(m_1,m_1,m_2,m_2,m_3,m_3)$,
is transformed to the new solution $(Y(\theta),X(\theta))^T$ in the variable true anomaly $\theta$
with $G=g=0$ with respect to the original Hamiltonian function $H$ of (\ref{new.H.function}), which is given by
\begin{equation}
Y(\theta)=\left(
\matrix{
\bar{Z}(\theta)\cr
\bar{W}(\theta)}
\right)
=\left(
\matrix{
0\cr
\sigma\cr
0\cr
0}
\right),
\quad
X(\theta)=\left(
\matrix{
\bar{z}(\theta)\cr
\bar{w}(\theta)}
\right)
=\left(
\matrix{
\sigma\cr
0\cr
0\cr
0}
\right).
\end{equation}

Moreover, the linearized Hamiltonian system at the elliptic Euler solution
${\xi}_0\equiv(Y(\theta),X(\theta))^T =$
\newline
$(0,\sigma,0,0,\sigma,0,0,0)^T\in\R^8$
depending on the true anomaly $\theta$ with respect to the Hamiltonian function
$H$ of (\ref{new.H.function}) is given by
\begin{equation}
\dot\zeta(\theta)=JB(\theta)\zeta(\theta),
\end{equation}
with
\begin{equation}
B(\theta)=H''(\theta,\bar{Z},\bar{W},\bar{z},\bar{w})|_{\bar\xi=\xi_0}
=\left(
\matrix{
I& O& -J& O\cr
O& I& O& -J\cr
J& O& H_{\bar{z}\bar{z}}(\theta,\xi_0)& O\cr
O& J& O& H_{\bar{w}\bar{w}}(\theta,\xi_0)
}
\right),
\end{equation}
and
\begin{equation}
H_{\bar{z}\bar{z}}(\theta,\xi_0)=\left(
\matrix{
-\frac{2-e\cos\theta}{1+e\cos\theta} & 0\cr
0 & 1
}
\right),
\quad
H_{\bar{w}\bar{w}}(\theta,\xi_0)=\left(
\matrix{
-\frac{2\delta-e\cos\theta}{1+e\cos\theta} & 0\cr
0 & \frac{\delta+e\cos\theta}{1+e\cos\theta}
}
\right),
\end{equation}
where
\begin{equation}\label{delta}
\delta=\frac{1}{\mu}\sum_{1\le i<j\le3}\frac{m_im_j(b_i-b_j)^2}{|a_{ix}-a_{jx}|^3}
=\frac{\sum_{1\le i<j\le3}\frac{m_im_j(b_i-b_j)^2}{|a_{ix}-a_{jx}|^3}}
{\sum_{1\le i<j\le3}\frac{m_im_j}{|a_{ix}-a_{jx}|}},
\end{equation}
and $H''$ is the Hession Matrix of $H$ with respect to its variable $\bar{Z}$,
$\bar{W}$, $\bar{z}$ and $\bar{w}$.
The corresponding quadratic Hamiltonian function is given by
\begin{eqnarray}
H_2(\theta,\bar{Z},\bar{W},\bar{z},\bar{w})
&=&\frac{1}{2}|\bar{Z}|^2+\bar{z}\cdot J\bar{Z}+\frac{1}{2}H_{\bar{z}\bar{z}}(\theta,\xi_0)|\bar{z}|^2
\nonumber\\
&&+\frac{1}{2}|\bar{W}|^2+\bar{w}\cdot J\bar{W}+\frac{1}{2}H_{\bar{w}\bar{w}}(\theta,\xi_0)|\bar{w}|^2.
\end{eqnarray}
\end{proposition}

{\bf Proof.} The proof is similar to those of Proposition 11.11 and Proposition 11.13 of \cite{Lon5}.
We just need to compute $H_{\bar{z}\bar{z}}(\theta,\xi_0)$, $H_{\bar{z}\bar{w}}(\theta,\xi_0)$
and $H_{\bar{w}\bar{w}}(\theta,\xi_0)$.

For simplicity, we omit all the upper bars on the variables of $H$ in (\ref{new.H.function}) in this proof.
By (\ref{new.H.function}), we have
\bea
H_z&=&JZ+\frac{p-r}{p}z-\frac{r}{\sigma}U_z(z,w),  \nn\\
H_w&=&JW+\frac{p-r}{p}w-\frac{r}{\sigma}U_w(z,w),  \nn\eea
and
\be\lb{Hessian}\left\{
\begin{array}{l}
H_{zz}=\frac{p-r}{p}I-\frac{r}{\sigma}U_{zz}(z,w),
\\
H_{zw}=H_{wz}=-\frac{r}{\sigma}U_{zw}(z,w),
\\
H_{ww}=\frac{p-r}{p}I-\frac{r}{\sigma}U_{ww}(z,w),
\end{array}\right. \ee
where we write $H_z$ and $H_{zw}$ etc to denote the derivative of $H$ with respect to $z$,
and the second derivative of $H$ with respect to $z$ and then $w$ respectively.
Note that all the items above are $2\times2$ matrices.

For $U_{ij}$ defined in (\ref{U_ij}) with $1\le i<j\le3$,
we have
\bea
\frac{\partial U_{ij}}{\partial z}(z,w) &=& -\frac{m_im_j(a_{ix}-a_{jx})}{|(a_{ix}-a_{jx})z+(b_i-b_j)w|^3}
\left[(a_{ix}-a_{jx})z+(b_i-b_j)w\right],  \nn\\
\frac{\partial U_{ij}}{\partial w}(z,w) &=& -\frac{m_im_j(b_i-b_j)}{|(a_{ix}-a_{jx})z+(b_i-b_j)w|^3}
\left[(a_{ix}-a_{jx})z+(b_i-b_j)w\right],  \nn\eea
and
\bea
&&\frac{\partial^2 U_{ij}}{\partial z^2}(z,w)=-\frac{m_im_j(a_{ix}-a_{jx})^2}{|(a_{ix}-a_{jx})z+(b_i-b_j)w|^3}I  \nn\\
&&\quad+3\frac{m_im_j(a_{ix}-a_{jx})^2}{|(a_{ix}-a_{jx})z+(b_i-b_j)w|^5}
\left[(a_{ix}-a_{jx})z+(b_i-b_j)w\right]\left[(a_{ix}-a_{jx})z+(b_i-b_j)w\right]^T,  \nn\\
&&\frac{\partial^2 U_{ij}}{\partial z\partial w}(z,w)=-\frac{m_im_j(a_{ix}-a_{jx})(b_i-b_j)}{|(a_{ix}-a_{jx})z+(b_i-b_j)w|^3}I \nn\\
&&\quad+3\frac{m_im_j(a_{ix}-a_{jx})(b_i-b_j)}{|(a_{ix}-a_{jx})z+(b_i-b_j)w|^5}
\left[(a_{ix}-a_{jx})z+(b_i-b_j)w\right]\left[(a_{ix}-a_{jx})z+(b_i-b_j)w\right]^T,  \nn\\
&&\frac{\partial^2 U_{ij}}{\partial w^2}(z,w)=-\frac{m_im_j(b_i-b_j)^2}{|(a_{ix}-a_{jx})z+(b_i-b_j)w|^3}I  \nn\\
&&\quad+3\frac{m_im_j(b_i-b_j)^2}{|(a_{ix}-a_{jx})z+(b_i-b_j)w|^5}
\left[(a_{ix}-a_{jx})z+(b_i-b_j)w\right]\left[(a_{ix}-a_{jx})z+(b_i-b_j)w\right]^T.  \nn\eea

Let
$$ K=\left(\matrix{2 & 0\cr
                   0 & -1}\right), \quad
K_1=\left(\matrix{1 & 0\cr
                  0 & 0}\right).  $$
Now evaluating these functions at the solution $\bar\xi_0=(0,\sigma,0,0,\sigma,0,0,0)^T\in\R^8$
 with $z=(\sigma,0)^T,w=(0,0)^T$, and summing them up,
we obtain
\begin{eqnarray}
\frac{\partial^2 U}{\partial z^2}\left|_{\xi_0}\right.&=&
\sum_{1\le i<j\le3}\frac{\partial^2 U_{ij}}{\partial z^2}\left|_{\xi_0}\right.
\nonumber\\
&=&\sum_{1\le i<j\le3}\left(-\frac{m_im_j(a_{ix}-a_{jx})^2}{|(a_{ix}-a_{jx})\sigma|^3}I
                        +3\frac{m_im_j(a_{ix}-a_{jx})^2}{|(a_{ix}-a_{jx})\sigma|^5}(a_{ix}-a_{jx})^2\sigma^2K_1\right)
\nonumber\\
&=&\frac{1}{\sigma^3}\left(\sum_{1\le i<j\le3}\frac{m_im_j}{|a_{ix}-a_{jx}|}\right)K
\nonumber\\
&=&\frac{\mu}{\sigma^3}K,  \label{U_zz}
\\
\frac{\partial^2 U}{\partial w^2}\left|_{\xi_0}\right.&=&
\sum_{1\le i<j\le3}\frac{\partial^2 U_{ij}}{\partial w^2}\left|_{\xi_0}\right.
\nonumber\\
&=&\sum_{1\le i<j\le3}\left(-\frac{m_im_j(b_i-b_j)^2}{|(a_{ix}-a_{jx})\sigma|^3}I
                        +3\frac{m_im_j(b_i-b_j)^2}{|(a_{ix}-a_{jx})\sigma|^5}(a_{ix}-a_{jx})^2\sigma^2K_1\right)
\nonumber\\
&=&\frac{1}{\sigma^3}\left(\sum_{1\le i<j\le3}\frac{m_im_j(b_i-b_j)^2}{|a_{ix}-a_{jx}|^3}\right)K
\nonumber\\
&=&\frac{\delta\mu}{\sigma^3}K, \label{U_ww}
\end{eqnarray}
where in the third equality of the first formula, we used (\ref{mu}),
and in the last equality of the second formula, we use the definition (\ref{delta}).
Similarly, we have
\begin{eqnarray}
\frac{\partial^2 U}{\partial z\partial w}\left|_{\xi_0}\right.&=&
\sum_{1\le i<j\le3}\frac{\partial^2 U_{ij}}{\partial z\partial w}\left|_{\xi_0}\right.
\nonumber\\
&=&\sum_{1\le i<j\le3}\left(-\frac{m_im_j(a_{ix}-a_{jx})(b_i-b_j)}{|(a_{ix}-a_{jx})\sigma|^3}I
                        +3\frac{m_im_j(a_{ix}-a_{jx})(b_i-b_j)}{|(a_{ix}-a_{jx})\sigma|^5}(a_{ix}-a_{jx})^2\sigma^2K_1\right)
\nonumber\\
&=&\left(\sum_{1\le i<j\le3}\frac{m_im_j(b_i-b_j)\cdot sign(a_{ix}-a_{jx})}{|(a_{ix}-a_{jx})|^2}\right)\frac{K}{\sigma^3}
\nonumber\\
&=&
\left(\frac{m_1m_2\sqrt{m_1m_2m_3}\alpha(\frac{1}{m_1}+\frac{1+x}{m_2})\cdot sign(-x\alpha)}{(-x\alpha)^2}
+\frac{m_2m_3[-\sqrt{m_1m_2m_3}\alpha(\frac{1+x}{m_2}+\frac{x}{m_3})]\cdot sign(-\alpha)}{(-\alpha)^2}\right.
\nonumber\\
&&+\left.\frac{m_3m_1\sqrt{m_1m_2m_3}\alpha(\frac{x}{m_3}-\frac{1}{m_1})\cdot sign((1+x)\alpha)}{((1+x)\alpha)^2}
\right)\frac{K}{\sigma^3}
\nonumber\\
&=&\frac{\sqrt{m_1m_2m_3}}{\alpha}\left(
-\frac{m_2+m_1+m_1x}{x^2}+(m_2+m_3)x+m_3+\frac{m_1x-m_3}{(1+x)^2}
\right)\frac{K}{\sigma^3}
\nonumber\\
&=&O,\label{U_zw}
\end{eqnarray}
where in the third equality, we used (\ref{a_ix}) and (\ref{b}), and in the last equality, we used
\begin{eqnarray}
&&-\frac{m_2+m_1x}{x^2}+(m_2+m_3)x+m_3+\frac{m_1x-m_3}{(1+x)^2}
\nonumber\\
&&=\frac{(m_2+m_3)x^5+(2m_2+3m_3)x^4+(m_2+3m_3)x^3-(3m_1+m_2)x^2-(3m_1+2m_2)x-(m_1+m_2)}{x^2(1+x)^2}
\nonumber\\
&&=0.
\end{eqnarray}

By (\ref{U_zz}), (\ref{U_ww}), (\ref{U_zw}) and (\ref{Hessian}), we have
\begin{eqnarray}
H_{zz}|_{\xi_0}&=&\frac{p-r}{p}I-\frac{r\mu}{\sigma^4}K=I-\frac{r}{p}I-\frac{r\mu}{p\mu}K
=I-\frac{r}{p}(I+K)
=\left(\matrix{-\frac{2-e\cos\theta}{1+e\cos\theta} & 0\cr
               0 & 1}\right),  \nn\\
H_{zw}|_{\xi_0}&=&-\frac{r}{\sigma}\frac{\partial^2U}{\partial z\partial w}|_{\xi_0}=O,  \nn\\
H_{ww}|_{\xi_0}&=&\frac{p-r}{p}I-\frac{r\delta\mu}{\sigma^4}K
=I-\frac{r}{p}I-\frac{r\delta\mu}{p\mu}K=I-\frac{r}{p}(I+\delta K)=
\left(\matrix{-\frac{2\delta-e\cos\theta}{1+e\cos\theta} & 0\cr
              0 & \frac{\delta+e\cos\theta}{1+e\cos\theta}}\right). \nn
\end{eqnarray}
Thus the proof is complete.\hb

We now want to obtain a simpler representation of $\delta$ of (\ref{delta}).
Plugging (\ref{a_ix}) and (\ref{b}) into (\ref{delta}), we have
\begin{eqnarray}\label{delta.simplify}
\delta&=&\frac{\frac{m_1m_2}{x^3}(\rho_1+\rho_2(1+x))^2
             +m_2m_3(\rho_2(1+x)+\rho_3 x)^2
             +\frac{m_3m_1}{(1+x)^3}(\rho_3 x-\rho_1)^2}
            {\frac{m_1m_2}{x}+m_2m_3+\frac{m_3m_1}{1+x}}
\nonumber\\
&=&\frac{m_3(1+x)^3(m_2+m_1+m_1x)^2+m_1x^3(1+x)^3(m_3+m_3x+m_2x)^2+m_2x^3(m_1x-m_3)^2}
        {x^2(1+x)^2[m_2m_3x^2+(m_1m_2+m_2m_3+m_3m_1)x+m_1m_2]}
\nonumber\\
&=&1+\frac{m_1(3x^2+3x+1)+m_3x^2(x^2+3x+3)}
          {x^2+m_2[(x+1)^2(x^2+1)-x^2]},
\end{eqnarray}
where $\rho_i$ are given by (\ref{rho_i}), and the last equality holds by Lemma 5.1 in the Appendix.
Note that $m_3=1-m_1-m_2$, the second term in the last equality of (\ref{delta.simplify})
is also defined in (18) of \cite{MaS1} (p. 317) and we use the same symbol $\beta$ of (\ref{1.4})
to denote it, and $\delta=\beta+1$. Then writing $H_{\bar{w}\bar{w}}(\theta,\bar\xi_0)$ in terms of
$\beta$ yields
\begin{equation}
H_{\bar{w}\bar{w}}(\theta,\xi_0)=\left(
\matrix{
-\frac{2\beta+2-e\cos\theta}{1+e\cos\theta} & 0\cr
0 & \frac{\beta+1+e\cos\theta}{1+e\cos\theta}
}
\right).
\end{equation}

Moreover, by the proof of Lemma 2 of \cite{MaS1}, we know that the full range of $\beta$ is $[0,7]$
when $m_1,m_2,m_3$ take all their possible values. Thus we have

\begin{proposition}\label{full.range}
The full range of the pair $(\beta,e)$ of the Euler elliptic orbit is the rectangle $[0,7]\times[0,1)$.
\end{proposition}

By Proposition \ref{linearized.Hamiltonian} , the essential part $\ga=\ga_{\bb,e}(t)$ of the
fundamental solution of the Euler orbit satisfies
\bea
\dot{\gamma}(t) &=& JB(t)\gamma(t),   \lb{2.17}\\
\gamma(0) &=& I_{4},    \lb{2.18}\eea
with
\be B(t)=\left(\matrix{1 & 0 & 0 & 1\cr
                       0 & 1 & -1 & 0 \cr
                       0 & -1 &\frac{-2\beta-2+e\cos(t)}{1+e\cos(t)} & 0 \cr
                       1 & 0 & 0 & \frac{\beta+1+e\cos(t)}{1+e\cos(t)} \cr}\right), \lb{2.19}\ee
where $e$ is the eccentricity, and $t$ is the truly anomaly.

Let
\be  J_2=\left(\matrix{ 0 & -1 \cr 1 & 0 \cr}\right), \qquad
    K_{\bb,e}(t)=\left(\matrix{\frac{2\beta+3}{1+e\cos(t)} & 0 \cr
                                     0 & -\frac{\beta}{1+e\cos(t)} \cr}\right),  \lb{2.20}\ee
and set
\be L(t,x,\dot{x})=\frac{1}{2}\|\dot{x}\|^2 + J_2x(t)\cdot\dot{x}(t) + \frac{1}{2}K_{\bb,e}(t)x(t)\cdot x(t),
       \qquad\quad  \forall\;x\in W^{1,2}(\R/2\pi\Z,\R^2),  \lb{2.21}\ee
where $a\cdot b$ denotes the inner product in $\R^2$. Obviously the origin in the configuration space is a
solution of the corresponding Euler-Lagrange system. By Legendrian transformation, the corresponding
Hamiltonian function is
$$   H(t,z)=\frac{1}{2}B(t)z\cdot z,\qquad \forall\; z\in\R^4.  $$

\subsection{A modification on the path $\ga_{\bb,e}(t)$}

In order to transform the Lagrangian system (\ref{2.19}) to a simpler linear operator corresponding to
a second order Hamiltonian system with the same linear stability as $\ga_{\bb,e}(2\pi)$, using $R(t)$
and $R_4(t)=\diag(R(t),R(t))$ as in Section 2.4 of \cite{HLS}, we let
\be  \xi_{\bb,e}(t) = R_4(t)\ga_{\bb,e}(t), \qquad \forall\; t\in [0,2\pi], (\bb,e)\in [0,7]\times [0,1). \lb{2.25}\ee
One can show by direct computation that
\be  \frac{d}{dt}\xi_{\bb,e}(t)
  = J \left(\matrix{I_2 & 0 \cr
                    0 & R(t)(I_2-K_{\bb,e}(t))R(t)^T \cr}\right)\xi_{\bb,e}(t). \lb{2.26}\ee
Note that $R_4(0)=R_4(2\pi)=I_4$, so $\ga_{\bb,e}(2\pi)=\xi_{\bb,e}(2\pi)$ holds. Then the linear stabilities
of the systems (\ref{2.18}) and (\ref{2.26}) are determined by the same matrix and thus is precisely the same.

By (\ref{2.25}) the symplectic paths $\ga_{\bb,e}$ and $\xi_{\bb,e}$ are homotopic to each other via the
homotopy $h(s,t)=R_4(st)\ga_{\bb,e}(t)$ for $(s,t)\in [0,1]\times [0,2\pi]$. Because $R_4(s)\ga_{\bb,e}(2\pi)$
for $s\in [0,1]$ is a loop in $\Sp(4)$ which is homotopic to the constant loop $\ga_{\bb,e}(2\pi)$,
$h(\cdot,2\pi)$ is contractible in $\Sp(4)$. Therefore by the proof of Lemma 5.2.2 on p.117 of \cite{Lon4},
the homotopy between $\ga_{\bb,e}$ and $\xi_{\bb,e}$ can be modified to fix the end point $\ga_{\bb,e}(2\pi)$
for all $s\in [0,1]$. Thus by the homotopy invariance of the Maslov-type index (cf. (i) of Theorem 6.2.7 on
p.147 of \cite{Lon4}) we obtain
\be  i_{\om}(\xi_{\bb,e}) = i_{\om}(\ga_{\bb,e}), \quad \nu_{\om}(\xi_{\bb,e}) = \nu_{\om}(\ga_{\bb,e}),
       \qquad \forall \,\omega\in\U, \; (\bb,e)\in [0,7]\times [0,1). \lb{2.27}\ee
Note that the first order linear Hamiltonian system (\ref{2.26}) corresponds to the following second order
linear Hamiltonian system
\be  \ddot{x}(t)=-x(t)+R(t)K_{\bb,e}(t)R(t)^Tx(t). \lb{2.28}\ee

For $(\bb,e)\in [0,7)\times [0,1)$, the second order differential operator corresponding to (\ref{2.28}) is given by
\bea  A(\bb,e)
&=& -\frac{d^2}{dt^2}I_2-I_2+R(t)K_{\bb,e}(t)R(t)^T  \nn\\
&=& -\frac{d^2}{dt^2}I_2-I_2+\frac{1}{2(1+e\cos t)}((3+\beta)I_2+3(1+\beta)S(t)),  \lb{2.29}\eea
where $S(t)=\left(\matrix{ \cos 2t & \sin 2t \cr
                           \sin 2t & -\cos 2t \cr}\right)$, defined on the domain $\ol{D}(\omega,2\pi)$
in (\ref{2.11}). Then it is self-adjoint and depends on the parameters $\bb$ and $e$. By Lemma
\ref{L2.3}, we have for any $\bb$ and $e$, the Morse index $\phi_{\om}(A(\bb,e))$ and nullity $\nu_{\om}(A(\bb,e))$
of the operator $A(\bb,e)$ on the domain $\ol{D}(\omega,2\pi)$ satisfy
\be  \phi_{\om}(A(\bb,e)) = i_{\om}(\xi_{\bb,e}), \quad \nu_{\om}(A(\bb,e)) = \nu_{\om}(\xi_{\bb,e}), \qquad
           \forall \,\om\in\U. \lb{2.30}\ee

In the rest of this paper, we shall use both of the paths $\ga_{\bb,e}$ and $\xi_{\bb,e}$ to study
the linear stability of $\ga_{\bb,e}(2\pi)=\xi_{\bb,e}(2\pi)$. Because of (\ref{2.27}), in many cases and
proofs below, we shall not distinguish these two paths.
Hence, if there is no confusion,
we will use $i_\om(\bb,e)$ and $\nu_\om(\bb,e)$ to represent $i_{\om}(\ga_{\bb,e})$
and $\nu_{\om}(\ga_{\bb,e})$ respectively.

\setcounter{equation}{0}
\section{Stability on the boundary of the unbounded rectangle $[0,\infty)\times [0,1)$}
\label{sec:3}

We start from the following lemma which will be used in sections 3 and 4. It is a special case
of Theorem 8.3.1 on p.188 of \cite{Lon4}, the details of whose proof is left to readers there
based on the ideas in the proofs of Theorems 8.2.1 and 8.2.2 on pp.184-185 of \cite{Lon4}. For
reader's conveniences, we give a detailed proof of this lemma here.

\begin{lemma}\label{Lm:Path.sum}
Let $\ga\in\mathcal{P}_{\tau}(4)$ satisfy
\be  \ga(\tau)\approx M_1\dm M_2   \lb{Path-1}\ee
with $M_1$, $M_2\in\Sp(2)$. Then there exist two paths $\ga_i\in\mathcal{P}_{\tau}(2)$ with
$\ga_i(\tau)=M_i$ for $i=1, 2$ such that we have
\be  i_1(\ga) \;=\; i_1(\ga_1) + i_1(\ga_2) \qquad{\it and}\qquad \ga\;\sim\;\ga_1\dm\ga_2.
          \lb{Path-2}\ee
\end{lemma}

{\bf Proof.} Firstly by Definition 5.2 below of $\ga(\tau)\approx M_1\dm M_2$ in (\ref{Path-1}),
there exists a continuous path $f\in C([0,\tau],\Om(\ga(\tau)))$ such that $f(0)=\ga(\tau)$ and
$f(\tau)=M_1\dm M_2$. We choose two paths $\xi$ and $\ga_2\in\mathcal{P}_{\tau}(2)$ satisfying
$\xi(\tau)=M_1$ and $\ga_2(\tau)=M_2$. Then $f\ast\ga(\tau) = M_1\dm M_2 = \xi\dm\ga_2(\tau)$.
Thus by Lemma 5.2.6 and Definition 5.2.7 on p.120 and Definition 5.4.2 on p.129 of \cite{Lon4},
there exists an integer $k\in\Z$ such that
$$  i_1(f\ast\ga) - (i_1(\xi) + i_1(\ga_2)) = 2k.  $$
Let $\phi_k(t) = R(2k\pi t/\tau)$ for $t\in [0,\tau]$. Define
$$  \ga_1(t) = \xi\ast\phi_k(t), \qquad \forall\;t\in [0,\tau]. $$
Then we obtain
$$  i_1(\ga_1) + i_1(\ga_2) = 2k + i_1(\xi) + i_1(\ga_2) = i_1(f\ast\ga).  $$
Thus by Theorem 6.2.4 on p.146 of \cite{Lon4} and the definition of the path $f$, we obtain
$$ \ga_1\dm\ga_2 \;\sim\; f\ast\ga \;\sim\; \ga,  $$
which completes the proof. \hb

By Proposition \ref{full.range}, we know the full range of $(\bb,e)$ is $[0,7]\times[0,1)$.
For convenience in the mathematical study, we extend the range of $(\bb,e)$ to $[0,\infty)\times [0,1)$.

Firstly, we need more precise information on indices and stabilities of $\ga_{\bb.e}$ at the boundary of the
$(\bb,e)$ rectangle $[0,\infty)\times [0,1)$.

\subsection{The boundary segment $\{0\}\times [0,1)$}

When $\beta=0$, this is the case if $m_1=0,\ m_2=1,\ m_3=0$, and the essential part of the fundamental solution of
Euler orbit is also the fundamental solution of the Keplerian orbits.
This is just the same case which has been discussed in Section 3.1 of \cite{HLS}.
We just cite the results here:
\begin{eqnarray}
i_{\om}(\ga_{0,e})=i_{\om}(\xi_{0,e})=
\left\{
\begin{array}{l}
0,\quad {\rm if}\;\; \om=1,\\
2,\quad {\rm if}\;\; \om \in \U\bs\{1\},
\end{array}
\right.
\lb{index.of.0e}
\\
\nu_{\om}(\ga_{0,e}) = \nu_{\om}(\xi_{0,e})
= \left\{
\begin{array}{l}
3, \quad {\rm if}\;\;\om = 1,
\\
0, \quad {\rm if}\;\;\om \in \U\bs\{1\}.
\end{array}
\right.
\lb{null.index.of.0e}
\end{eqnarray}

\subsection{The boundary $[0,\infty)\times \{0\}$}

In this case $e=0$. It is considered in (A) of Subsection 3.1 of \cite{HLS} when $\bb=0$.
Below, we shall first recall the
properties of eigenvalues of $\ga_{\bb,0}(2\pi)$. Then we carry out the computations of normal
forms of $\ga_{\bb,0}(2\pi)$, and $\pm 1$ indices $i_{\pm 1}(\ga_{\bb,0})$ of the path
$\ga_{\bb,0}$ for all $\bb\in [0,\infty)$, which are new.

In this case, the essential part of the motion (\ref{2.17})-(\ref{2.19}) becomes an ODE
system with constant coefficients:
\be B = B(t) = \left(\matrix{1 & 0 & 0 & 1\cr
                             0 & 1 & -1 & 0 \cr
                             0 & -1 & -2\bb-2 & 0 \cr
                             1 & 0 & 0 & \bb+1 \cr}\right).  \lb{3.10}\ee
The characteristic polynomial $\det(JB-\lambda I)$ of $JB$ is given by
\be \lambda^4 + (1-\bb)\lambda^2-\bb(2\bb+3) = 0.  \lb{3.11}\ee
Letting $\aa=\lambda^2$, the two roots of the quadratic polynomial $\aa^2 + (1-\bb)\aa -\bb(2\bb+3)$
are given by $\aa_1=\frac{\bb-1+\sqrt{9\bb^2+10\bb+1}}{2}\ge0$
and $\aa_2=\frac{\bb-1-\sqrt{9\bb^2+10\bb+1}}{2}<0$.
Therefore the four roots of the polynomial
(\ref{3.11}) are given by
\begin{eqnarray}
\aa_{1,\pm} &=&\pm\sqrt{\aa_1}= \pm\sqrt{\frac{\bb-1+\sqrt{9\bb^2+10\bb+1}}{2}}\in\R,\lb{3.12a}
\\
\aa_{2,\pm} &=&\pm\sqrt{-1}\sqrt{-\aa_2}= \pm\sqrt{-1}\sqrt{\frac{-\bb+1+\sqrt{9\bb^2+10\bb+1}}{2}}. \lb{3.12}
\end{eqnarray}
Moreover, when $\bb\ge0$, we have
\begin{eqnarray}
\frac{d\aa_1}{d\bb}=\frac{1}{2}+\frac{9\bb+5}{2\sqrt{9\bb^2+10\bb+1}}>0,\label{da1.db}\\
\frac{d\aa_2}{d\bb}=\frac{1}{2}-\frac{9\bb+5}{2\sqrt{9\bb^2+10\bb+1}}<0.\label{da2.db}
\end{eqnarray}

{\bf (A)  Eigenvalues of $\ga_{\bb,0}(2\pi)$ for $\bb\in [0,\infty)$.}

When $\bb\ge 0$, by (\ref{3.12a}) and (\ref{3.12}), we get the four characteristic multipliers
of the matrix $\ga_{\bb,0}(2\pi)$
\begin{eqnarray}
\rho_{1,\pm}(\beta) = e^{2\pi\aa_{1,\pm}} =  e^{\pm 2\pi\sqrt{\aa_1}}\in\R^+, \lb{3.13a}
\\
\rho_{2,\pm}(\beta) = e^{2\pi\aa_{2,\pm}} =  e^{\pm 2\pi\sqrt{-1}\th(\bb)}, \lb{3.13}
\end{eqnarray}
where
\be  \th(\bb) = \sqrt{\frac{-\bb+1+\sqrt{9\bb^2+10\bb+1}}{2}}. \lb{3.14}\ee
By (\ref{da2.db}) and (\ref{3.14}), we know that $\th(\bb)$ is increasing with respect to $\bb$ when $\bb\ge0$.

From (\ref{3.14}), $\th(0)=1$.
Then for any $\th\ge1$, we denote by $\bb_\th\ge0$ the $\bb$ value satisfying $\th(\bb)=\th$,
and we obtain
$$
\th=\sqrt{\frac{-\bb_\th+1+\sqrt{9\bb_\th^2+10\bb_\th+1}}{2}},
$$
and hence
\begin{equation}
\bb_\th=\frac{\th^2-3+\sqrt{9\th^4-14\th^2+9}}{4},\quad\quad \th\ge1. \label{om_th}
\end{equation}
Moreover, when $\th\ge1$, we have
\begin{equation}
\frac{d\bb_\th}{d\th}=\frac{2\th+\frac{2\th(9\th^2-7)}{\sqrt{9\th^4-14\th^2+9}}}{4}>0. \label{dom_th.dth}
\end{equation}
For later use, we write $\bb_\th$ for $\th=n$ and $\th=n+\frac{1}{2}$, $n\in\N$ as
\begin{equation}
\hat\bb_n=\frac{n^2-3+\sqrt{9n^4-14n^2+9}}{4},\quad\quad n=1,\;2,\;3...\label{om_n}
\end{equation}
and
\begin{equation}
\hat\bb_{n+\frac{1}{2}}=\frac{(n+\frac{1}{2})^2-3+\sqrt{9(n+\frac{1}{2})^4-14(n+\frac{1}{2})^2+9}}{4},
    \quad\quad n=1,\;2,\;3...\label{om_n.5}
\end{equation}
where we have used the symbol hat to denote these special values of $\bb$.
Moreover, from (\ref{om_n}) we have
\begin{eqnarray}\label{bb_n.approx}
\hat\bb_n&=&\frac{n^2-3+\sqrt{9n^4-14n^2+9}}{4}\nonumber\\
&=&n^2-\frac{3n^2+3-\sqrt{9n^4-14n^2+9}}{4}\nonumber\\
&=&n^2-\frac{32n^2}{4(3n^2+3+\sqrt{9n^4-14n^2+9})}\nonumber\\
&=&n^2-\frac{8}{3+\frac{3}{n^2}+\sqrt{9-\frac{14}{n^2}+\frac{9}{n^4}}}\nonumber\\
&\approx&n^2-\frac{4}{3},
\end{eqnarray}
when $n$ is large enough.
By (\ref{dom_th.dth}), we have
\begin{equation}\label{rank}
0=\hat\bb_1<\hat\bb_{\frac{3}{2}}<\hat\bb_2<\hat\bb_{\frac{5}{2}}<...<\hat\bb_n<\hat\bb_{n+\frac{1}{2}}<...
\end{equation}

Specially, we obtain the following results:

(i) When $\bb=\hat\bb_1=0$, we have $\sg(\ga_{0,0}(2\pi)) = \{1, 1, 1, 1\}$.

When $\bb>0$, by (\ref{da1.db}), (\ref{da2.db}) and (\ref{3.13a}),
we have $\aa_1>0$, and hence $\rho_{1,\pm}(\bb)=e^{\pm2\pi \sqrt{\aa_1}}\subset \R\bs\U$.

(ii) Let $i\in\N$.
When $\hat\bb_i<\bb<\hat\bb_{i+\frac{1}{2}}$, the angle $\th(\bb)$ in (\ref{3.14}) increases strictly from $i$ to $i+\frac{1}{2}$ as $\bb$ increases
from $\hat\bb_i$ to $\hat\bb_{i+\frac{1}{2}}$.
Therefore $\rho_{2,+}(\bb)=e^{2\pi \sqrt{-1}\th(\bb)}$ runs from $1$ to $-1$ counterclockwise along
the upper semi-unit circle in the complex plane $\C$ as $\bb$ increases from $\hat\bb_i$ to $\hat\bb_{i+\frac{1}{2}}$.
Correspondingly
$\rho_{2,-}(\bb)=e^{-2\pi \sqrt{-1}\th(\bb)}$ runs from $1$ to $-1$ clockwise along the lower semi-unit circle in
$\C$ as $\bb$ increases from $\hat\bb_i$ to $\hat\bb_{i+\frac{1}{2}}$.
Thus specially we obtain
$\rho_{2,\pm}(\bb)\subset \U\bs\R$ for all $\bb\in (\hat\bb_i,\hat\bb_{i+\frac{1}{2}})$.

(iii) When $\bb=\hat\bb_{i+\frac{1}{2}}$, we have $\th(\hat\bb_{i+\frac{1}{2}})=i+\frac{1}{2}$. Therefore we obtain
$\rho_{2,\pm}(\hat\bb_{i+\frac{1}{2}})=e^{\pm \sqrt{-1} \pi} = -1$.

(iv) When $\hat\bb_{i+\frac{1}{2}}<\bb<\hat\bb_{i+1}$, the angle $\th(\bb)$ increases strictly from $i+\frac{1}{2}$ to $i+1$ as $\bb$ increase
from $\hat\bb_{i+\frac{1}{2}}$ to $\hat\bb_{i+1}$. Thus $\rho_{2,+}(\bb)=e^{2\pi \sqrt{-1}\th(\bb)}$ runs from $-1$ to $1$
counterclockwise along the lower semi-unit circle in $\C$ as $\bb$ increases from $\hat\bb_{i+\frac{1}{2}}$ to $\hat\bb_{i+1}$. Correspondingly
$\rho_{2,-}(\bb)=e^{-2\pi \sqrt{-1}\th(\bb)}$ runs from $-1$ to $1$ clockwise along the
upper semi-unit circle in $\C$ as $\bb$ increases from $\hat\bb_{i+\frac{1}{2}}$ to $\hat\bb_{i+1}$.
Thus we obtain
$\rho_{2,\pm}(\bb) \subset \U\bs\R$ for all $\bb\in (\hat\bb_{i+\frac{1}{2}},\hat\bb_{i+1})$.

(v) When $\bb=\hat\bb_{i+1}$, we obtain $\th(\hat\bb_{i+1})=i+1$, and then we have double eigenvalues
$\rho_{2,\pm}(\hat\bb_{i+1}) = 1$.

{\bf (B)  Indices $i_1(\ga_{\bb,0})$ of $\ga_{\bb,0}(2\pi)$ for $\bb\in [0,\infty)$.}

Define
\begin{equation}\lb{f0}
f_{0,1}=R(t)\left(\matrix{1\cr 0}\right),\quad
f_{0,2}=R(t)\left(\matrix{0\cr 1}\right),
\end{equation}
and
\begin{equation}\lb{fn}
f_{n,1}=R(t)\left(\matrix{\cos nt\cr 0}\right),\quad
f_{n,2}=R(t)\left(\matrix{0\cr \cos nt}\right),\quad
f_{n,3}=R(t)\left(\matrix{\sin nt\cr 0}\right),\quad
f_{n,4}=R(t)\left(\matrix{0\cr \sin nt}\right),
\end{equation}
for $n\in\N$.
Then $f_{0,1},\;f_{0,2}$ and $f_{n,1},\;f_{n,2}\;f_{n,3},\;f_{n,4}\;n\in\N$
form an orthogonal basis of $\overline{D}(1,2\pi)$.
By (\ref{2.29}) and $\frac{dR(t)}{dt}=JR(t)$, computing $A(\bb,0)f_{n,1}$ yields
\begin{eqnarray}
A(\bb,0)f_{n,1}&=&[-\frac{d^2}{dt^2}I_2-I_2+R(t)K_{0,e}(t)R(t)^T]R(t)\left(\matrix{\cos nt\cr 0}\right)
\nonumber\\
&=&R(t)\left(\matrix{(n^2+2\bb+3)\cos nt\cr 2n\sin nt}\right)
\nonumber\\
&=&(n^2+2\bb+3)f_{n,1}+2nf_{n,4}.
\end{eqnarray}
Similarly, we have
\begin{eqnarray}
\left(\matrix{A(\bb,0) & O\cr O & A(\bb,0)}\right)\left(\matrix{f_{0,1}\cr f_{0,2}}\right)&=&
\left(\matrix{2\bb+3 & 0\cr 0 & -\bb}\right)
\left(\matrix{f_{0,1}\cr f_{0,2}}\right),\lb{A1}
\\
\left(\matrix{A(\bb,0) & O\cr O & A(\bb,0)}\right)\left(\matrix{f_{n,1}\cr f_{n,4}}\right)&=&
\left(\matrix{n^2+2\bb+3 & 2n\cr 2n & n^2-\bb}\right)
\left(\matrix{f_{n,1}\cr f_{n,4}}\right),\lb{A2}
\\
\left(\matrix{A(\bb,0) & O\cr O & A(\bb,0)}\right)\left(\matrix{f_{n,3}\cr f_{n,2}}\right)&=&
\left(\matrix{n^2+2\bb+3 & -2n\cr -2n & n^2-\bb}\right)
\left(\matrix{f_{n,3}\cr f_{n,2}}\right),\lb{A3}
\end{eqnarray}
for $n\in\N$.
Denote
\begin{equation}
B_0=\left(\matrix{2\bb+3 & 0\cr 0 & -\bb}\right),\quad
B_n=\left(\matrix{n^2+2\bb+3 & 2n\cr 2n & n^2-\bb}\right),\quad
\tilde{B}_n=\left(\matrix{n^2+2\bb+3 & -2n\cr -2n & n^2-\bb}\right).
\end{equation}
Denote the characteristic polynomial of $B_n$ and $\tilde{B}_n$ by $p_n(\lambda)$ and $\tilde{p}_n(\lambda)$
respectively, then we have
\begin{equation}
p_n(\lambda)=\tilde{p}_n(\lambda)=\lambda^2-(2n^2+\bb+3)\lambda-[2\bb^2-(n^2-3)\bb-n^2(n^2-1)]
\end{equation}

Let $i\in\N\;i>1$, fix $\bb=\hat\bb_i$, then $p_n(0)=\tilde{p}_n(0)=0$ iff $n=i$.
Moreover, we have $p_n(0)=\tilde{p}_n(0)<0$ if $n<i$,
and $p_n(0)=\tilde{p}_n(0)>0$ if $n>i$.
Thus both $B_i$ and $\tilde{B}_i$ have a zero and a positive eigenvalues;
both $B_n$ and $\tilde{B}_n$ with $n<i$ have a negative and a positive eigenvalues;
both $B_n$ and $\tilde{B}_n$ with $n>i$ have two positive eigenvalues.
Notice that $B_0$ has a negative and a positive eigenvalues.
Then we have $i_1(\ga_{\hat\bb_i,0})=2i-1$ and $\nu_1(\ga_{\hat\bb_i,0})=2$.

When $\hat\bb_i<\bb<\hat\bb_{i+1}$, then $p_n(0)=\tilde{p}_n(0)\ne0$.
Similarly to the above argument, we have $p_n(0)=\tilde{p}_n(0)<0$ if $n\le i$,
and $p_n(0)=\tilde{p}_n(0)>0$ if $n>i$.
Thus both $B_n$ and $\tilde{B}_n$ with $n\le i$ have a negative and a positive eigenvalues;
both $B_n$ and $\tilde{B}_n$ with $n>i$ have two positive eigenvalues.
Notice that $B_0$ has a negative and a positive eigenvalues,
we have $i_1(\ga_{\bb,0})=2i+1$ and $\nu_1(\ga_{\bb,0})=0$.

Therefore, we have
\begin{eqnarray}
&& i_1(\ga_{\bb,0}) = \left\{\matrix{  
                 0, &  {\rm if}\;\;\bb=\hat\bb_1=0, \cr
                 3, &  {\rm if}\;\;\bb\in(\hat\bb_1,\hat\bb_2], \cr
                 ...,\cr
                 2n+1, &  {\rm if}\;\;\bb\in(\hat\bb_n,\hat\bb_{n+1}], \cr
                 ...\cr}\right.\lb{1-index.of.b0}
\\
&& \nu_1(\ga_{\bb,0}) = \left\{\matrix{
                 3, &  {\rm if}\;\;\bb=\hat\bb_1=0, \cr
                 2, &  {\rm if}\;\;\bb=\hat\bb_n,\;n\ge2, \cr
                 0, &  {\rm if}\;\;\bb\ne\hat\bb_1,\hat\bb_2,...\hat\bb_n,... \cr}\right. \lb{null.1-index.of.b0}
\end{eqnarray}
where the case of $\bb=\hat\bb_1=0$ follows from (\ref{index.of.0e}) and (\ref{null.index.of.0e}).

{\bf (C)  Indices $i_{\om}(\ga_{\bb,0}),\;\om\ne1$ for $\bb\in [0,\infty)$.}

By a similar arguments in (B), we can compute the eigenvalues of $A(\bb,0)$ in the domain
$\overline{D}(-1,2\pi)$, and hence the $-1$-indices of $\ga_{\bb,0}$. Especially, when
$\bb=\hat\bb_{n+{1\over2}}$, $A(\bb,0)$ has eigenvalue $-1$ with multiplicity $2$. Thus
\be\label{i.-1}
i_{-1}(\ga_{\hat\bb_{n+1\slash2,0}}(2\pi))=2.
\ee


From the above discussions, when $\bb\ge 0$, by (\ref{3.12a})-(\ref{da2.db}) and (i)-(v) in Part
(A), $\ga_{\bb,0}(2\pi)$ possesses one pair of positive hyperbolic characteristic multipliers
$\rho_{1,\pm}(\bb)$ given by (\ref{3.13a}), and one pair of elliptic characteristic multipliers
$\rho_{2,\pm}(\bb)$ on the unit circle given by (\ref{3.13}). Therefore by Theorem 1.7.3 on p.36
of \cite{Lon4}, we have
\be  \ga_{\bb,0}(2\pi)\approx D(e^{2\pi\sqrt{\aa_1(\bb)}})\dm M(2\pi)  \lb{3.A1}\ee
for some matrix $M(2\pi)\in\Sp(2)$ satisfying
\be  M(2\pi)
= \left\{
\begin{array}{l}
I_2,\quad {\rm if}\;\bb=\hat\bb_n,n\in\N,\\
-I_2,\quad {\rm if}\; \bb=\hat\bb_{n+\frac{1}{2}},n\in\N,\\
R(2\pi\th(\bb))\;{\rm or}\;R(-2\pi\th(\bb)),\quad {\rm if}\; \bb\ne\hat\bb_n,\;\hat\bb_{n+\frac{1}{2}},\;\forall n\in\N,
\end{array}
\right. \lb{3.A2}
\end{equation}
where we have used (i)-(v) in Part (A) again.

By Lemma \ref{Lm:Path.sum} there exists a path $M\in\mathcal{P}_{2\pi}(2)$ connecting $M(0)=I_2$ to $M(2\pi)$
such that the path $\ga_{\bb,0}(t)$ is homotopic to the path $D(e^{t\sqrt{\aa_1(\bb)}})\diamond M(t)$ defined
for $t\in [0,2\pi]$.

By the properties of splitting numbers in Chapter 9 of \cite{Lon4}, for $\hat\bb_n<\bb<\hat\bb_{n+\frac{1}{2}}$
and $\om=-1$, we obtain
\begin{eqnarray} i_{-1}(\ga_{\bb,0})
&=& i_{1}(\ga_{\bb,0})+S_{\ga_{\bb,0}(2\pi)}^+(1)-S_{\ga_{\bb,0}(2\pi)}^-(e^{\sqrt{-1}2\pi[\th(\bb)-n]})
    +S_{\ga_{\bb,0}(2\pi)}^+(e^{\sqrt{-1}2\pi[\th(\bb)-n]})-S_{\ga_{\bb,0}(2\pi)}^-(-1)  \nn\\
&=& i_{1}(\ga_{\bb,0})-S_{\ga_{\bb,0}(2\pi)}^-(e^{\sqrt{-1}2\pi[\th(\bb)-n]})
    +S_{\ga_{\bb,0}(2\pi)}^+(e^{\sqrt{-1}2\pi[\th(\bb)-n]})   \nn\\
&=& i_{1}(\ga_{\bb,0})-S_{D(exp(2\pi\sqrt{\aa_1}))}^-(e^{\sqrt{-1}2\pi[\th(\bb)-n]})
    - S_{M(2\pi)}^-(e^{\sqrt{-1}2\pi[\th(\bb)-n]})   \nn\\
& & +S_{D(exp(2\pi\sqrt{\aa_1}))}^+(e^{\sqrt{-1}2\pi[\th(\bb)-n]})+S_{M(2\pi)}^+(e^{\sqrt{-1}2\pi[\th(\bb)-n]}) \nn\\
&=& i_{1}(\ga_{\bb,0})-S_{M(2\pi)}^-(e^{\sqrt{-1}2\pi[\th(\bb)-n]})+S_{M(2\pi)}^+(e^{\sqrt{-1}2\pi[\th(\bb)-n]}) \nn\\
&=&\left\{
\begin{array}{l}
i_{1}(\ga_{\bb,0})-1=2n,\quad{\rm if}\;M(2\pi)=R(2\pi\th(\bb)), \\
i_{1}(\ga_{\bb,0})+1=2n+2,\quad{\rm if}\;M(2\pi)=R(-2\pi\th(\bb)), \end{array} \right.  \lb{3.A3}
\end{eqnarray}
where the first equality follows from (\ref{5.A1}) below, the second equality follows from the fact
$1\not\in \sg(\ga_{\bb,0}(2\pi))$ by (\ref{3.A1}) and the third case of (\ref{3.A2}), the third equality
follows from (\ref{3.A1}), the forth equality follows from the fact
$\U\cap\sg(D(exp(2\pi\sqrt{\aa_1})))=\emptyset$, and in the last step we have used
(\ref{1-index.of.b0})-(\ref{null.1-index.of.b0}).

Similarly, when $\hat\bb_{n+\frac{1}{2}}<\bb<\hat\bb_{n+1}$, we have
\begin{eqnarray} i_{-1}(\ga_{\bb,0})
&=& i_{1}(\ga_{\bb,0})+S_{\ga_{\bb,0}(2\pi)}^+(1)-S_{\ga_{\bb,0}(2\pi)}^-(e^{\sqrt{-1}2\pi[n+1-\th(\bb)]})
    + S_{\ga_{\bb,0}(2\pi)}^+(e^{\sqrt{-1}2\pi[n+1-\th(\bb)]})-S_{\ga_{\bb,0}(2\pi)}^-(-1)   \nn\\
&=& i_{1}(\ga_{\bb,0})-S_{M(2\pi)}^-(e^{\sqrt{-1}2\pi[n+1-\th(\bb)]})
    + S_{M(2\pi)}^+(e^{\sqrt{-1}2\pi[n+1-\th(\bb)]})    \nn\\
&=&\left\{
\begin{array}{l}
i_{1}(\ga_{\bb,0})+1=2n+2,\quad{\rm if}\;M(2\pi)=R(2\pi\th(\bb)), \\
i_{1}(\ga_{\bb,0})-1=2n,\quad{\rm if}\;M(2\pi)=R(-2\pi\th(\bb)). \end{array}\right.  \label{i.-1.computation}
\end{eqnarray}

If $M(2\pi)=R(-2\pi\th(\bb))$ for $\hat\bb_n<\bb<\hat\bb_{n+\frac{1}{2}}$,
we have $i_{-1}(\ga_{\bb,0})=2n+2$.
By (\ref{i.-1}) and the non-decreasing of $i_{-1}(\ga_{\bb,0})$ with respect to $\bb$ of Lemma 4.2 below,
we must have
$i_{-1}(\ga_{\hat\bb_{n+1\slash2}+\epsilon,0})
= i_{-1}(\ga_{\hat\bb_{n+1\slash2},0})+\nu_{-1}(\ga_{\hat\bb_{n+1\slash2},0}) \ge 2n+4$,
which contradicts (\ref{i.-1.computation}). Similarly, we cannot have $M(2\pi)=R(-2\pi\th(\bb))$ for
$\hat\bb_{n+\frac{1}{2}}<\bb<\hat\bb_{n+1}$, too. Thus we must have $M(2\pi)=R(2\pi\th(\bb))$ when
$\bb\ne\hat\bb_n,\;\hat\bb_{n+\frac{1}{2}},\;\forall n\in\N$.

Therefore,
\begin{eqnarray}
&& i_{-1}(\ga_{\bb,0}) = \left\{\matrix{
                 2, &  {\rm if}\;\;\bb\in[0,\hat\bb_{\frac{3}{2}}], \cr
                 4, &  {\rm if}\;\;\bb\in(\hat\bb_{\frac{3}{2}},\hat\bb_{\frac{5}{2}}], \cr
                 ...,\cr
                 2n, &  {\rm if}\;\;\bb\in(\hat\bb_{n-\frac{1}{2}},\hat\bb_{n+\frac{1}{2}}],\;n\ge2, \cr
                 ...\cr}\right. \lb{3.A4}\\
&& \nu_{-1}(\ga_{\bb,0}) = \left\{\matrix{
                 2, &  {\rm if}\;\;\bb=\hat\bb_{n+\frac{1}{2}},\;n\in\N, \cr
                 0, &  {\rm if}\;\;\bb\ne\hat\bb_{\frac{3}{2}},\hat\bb_{\frac{5}{2}},...\bb_{n+\frac{1}{2}},... \cr}, \right.
\lb{null.-1.index.of.b0}
\end{eqnarray}
where in (\ref{3.A4}) we have used the left continuity of the index functions at the degenerate
points to get their values at $\bb=\hat{\bb}_n$ or $\hat{\bb}_{n+1/2}$ (cf. Definition 5.4.2 on
p.129 of \cite{Lon4}).

For any real number $\th_0$ such that $0<\th_0<\frac{1}{2}$. Let $\om_0=e^{2\pi\th_0\sqrt{-1}}$,

Similarly, for $\omega\in{\bf U}\backslash\{1,-1\}$, $i_{\om_0}(\ga_{\bb,0})$ can be computed
using the decreasing property of the index proved in Corollary \ref{C4.5}.

\setcounter{equation}{0}
\section{The degeneracy curves of elliptic Euler solutions}\label{sec:4}

\subsection{The increasing of $\om$-indeces of elliptic Euler solutions}

For convenience, we define
\bea
A_1(e)&=&-\frac{d^2}{dt^2}-1+\frac{1}{1+e\cos t},
\\
A(-1,e)&=&-\frac{d^2}{dt^2}I_2-I_2+\frac{I_2}{1+e\cos t}=A_1(e)\oplus A_1(e).\label{A.-1.e}
\eea

For $(\bb,e)\in [0,\infty)\times [0,1)$, let $\bar{A}(\bb,e)=\frac{A(\bb,e)}{\bb+1}$. Using
(\ref{2.29}) we can rewrite $A(\bb,e)$ as follows
\begin{eqnarray} A(\bb,e) &=& -\frac{d^2}{dt^2}I_2-I_2+\frac{I_2}{1+e\cos t}+\frac{\bb+1}{2(1+e\cos t)}(I_2+3S(t))
\nonumber\\
&=& (\beta+1)\bar{A}(\bb,e), \lb{4.11}
\end{eqnarray}
where we define
\be\label{bar.A}
\bar{A}(\bb,e)=\frac{-\frac{d^2}{dt^2}I_2-I_2+\frac{I_2}{1+e\cos t}}{\bb+1}+\frac{I_2+3S(t)}{2(1+e\cos t)}
=\frac{A(-1,e)}{\bb+1}+\frac{I_2+3S(t)}{2(1+e\cos t)}.
\ee
Therefore we have
\bea
\phi_\omega(A(\bb,e)) &=& \phi_\omega(\bar{A}(\bb,e)), \lb{4.12}\\
\nu_\omega(A(\bb,e)) &=& \nu_\omega(\bar{A}(\bb,e)).   \lb{4.13}\eea

In \cite{HO}, Hu and Ou proved that the operator $-\frac{d^2}{dt^2}-1+\frac{\bb}{1+e\cos t}$ is positive definite
for $\bb>1$. Moreover, we have
\begin{lemma}
For $0\le e<1$, there holds

(i)$A_1(e)$ and $A(-1,e)$ are non-negative definite for the $\omega=1$
boundary condition, and
\bea
\ker{A_1(e)}&=&\{c(1+e\cos t)|c\in\C\},\label{Ker.A1}
\\
\ker{A(-1,e)}&=&\left\{\left(\matrix{c_1(1+e\cos t)\cr c_2(1+e\cos t)}\right)\bigg|c_1,c_2\in\C\right\},\label{Ker.A.-1.e}
\eea

(ii)$A_1(e)$ and $A(-1,e)$ are positive definite for any $\omega\ne1$ boundary condition.
\end{lemma}
{\bf Proof.} By (\ref{A.-1.e}), we just need to prove the results for $A_1(e)$.
Let $x(t)\not\equiv0\in\overline{D}(\omega, 2\pi)$, then
\be
y(t)=\frac{x(t)}{1+e\cos t}\in\overline{D}(\omega, 2\pi).
\ee
Then we have
\bea
\<A_1(e)x(t),x(t)\>&=&\int_0^{2\pi}[|x'(t)|^2-\frac{e\cos t}{1+e\cos t}|x(t)|^2]dt
\nonumber\\
&=&\int_0^{2\pi}[(1+e\cos t)^2|y'(t)|^2+e^2\sin^2t|y(t)|^2
   -e\sin t(1+e\cos t)(y(t)\overline{y'(t)}+\overline{y(t)}y'(t))]dt
\nonumber\\
&&-\int_0^{2\pi}{e\cos t}(1+e\cos t)|y(t)|^2dt
\nonumber\\
&=&\int_0^{2\pi}[(1+e\cos t)^2|y'(t)|^2+e^2\sin^2t|y(t)|^2
   -e\sin t(1+e\cos t)(y(t)\overline{y'(t)}+\overline{y(t)}y'(t))]dt
\nonumber\\
&&-\int_0^{2\pi}e(1+e\cos t)|y(t)|^2d\sin t
\nonumber\\
&=&\int_0^{2\pi}[(1+e\cos t)^2|y'(t)|^2+e^2\sin^2t|y(t)|^2
  -e\sin t(1+e\cos t)(y(t)\overline{y'(t)}+\overline{y(t)}y'(t))]dt
\nonumber\\
&&+\int_0^{2\pi}\sin td[e(1+e\cos t)|y(t)|^2]
\nonumber\\
&=&\int_0^{2\pi}[(1+e\cos t)^2|y'(t)|^2+e^2\sin^2t|y(t)|^2
   -e\sin t(1+e\cos t)(y(t)\overline{y'(t)}+\overline{y(t)}y'(t))]dt
\nonumber\\
&&+\int_0^{2\pi}[-e^2\sin^2t|y(t)|^2+e\sin t(1+e\cos t)(y(t)\overline{y'(t)}+\overline{y(t)}y'(t))]dt
\nonumber\\
&=&\int_0^{2\pi}(1+e\cos t)^2|y'(t)|^2dt
\nonumber\\
&\ge&0,
\eea
where the last equality holds if and only if $y(t)\equiv c$ for some constant $c\ne0$.
In such case, we have $x(0)=x(2\pi)=c\ne0$, which can be happen when $\omega=1$ but not for $\omega\in\U\backslash{1}$.
Therefore, $A_1(e)$ is positive definite for any $\omega\ne1$
boundary condition;
non-negative definite for the $\omega=1$ boundary condition, and in such case, (\ref{Ker.A1}) holds.
\hb

Now motivated by Lemma 4.4 in \cite{HLS} and modifying its proof to the Euler case,
we get the following important lemma:

\begin{lemma}\label{Lemma:increasing.of.index}
(i) For each fixed $e\in [0,1)$, the operator $\bar{A}(\bb,e)$ is non-increasing
with respect to $\beta\in [0,+\infty)$ for any fixed $\omega\in\U$. Specially
\be  \frac{\pt}{\pt\beta}\bar{A}(\beta,e)|_{\bb=\bb_0} = -\frac{1}{(\bb_0+1)^2}A(-1,e),  \lb{4.14}\ee
is a non-negative definite operator for $\bb_0\in [0,\infty)$.

(ii) For every eigenvalue $\lm_{\bb_0}=0$ of $\bar{A}(\bb_0,e_0)$ with $\om\in\U$ for some
$(\bb_0,e_0)\in [0,\infty)\times [0,1)$, there holds
\be \frac{d}{d\bb}\lm_{\bb}|_{\bb=\bb_0} < 0.  \lb{4.15}\ee

(iii)  For every $(\bb,e)\in(0,\infty)\times[0,1)$ and $\om\in\U$,
there exist $\epsilon_0=\epsilon_0(\bb,e)>0$ small enough such that for all $\epsilon\in(0,\epsilon_0)$ there holds
\be
i_\om(\ga_{\bb+\epsilon,e})-i_\om(\ga_{\bb,e})=\nu_\om(\ga_{\bb,e}).
\ee
\end{lemma}

{\bf Proof.} If we have (\ref{4.15}), (iii) can be proved by using the same techniques in the proof
of the first part of Proposition 6.1 in \cite{HLS}.
So it suffices to prove (ii).
Let $x_0=x_0(t)$ with unit norm such that
\be  \bar{A}(\bb_0,e_0)x_0=0.     \lb{4.16}\ee
Fix $e_0$, then $\bar{A}(\bb,e_0)$ is an analytic path of non-increasing self-adjoint operators
with respect to $\beta$. Following Kato (\cite{Ka}, p.120 and p.386), we can choose a smooth path of
unit norm eigenvectors $x_{\bb}$ with $x_{\beta_0}= x_0$ belonging to a smooth path of real eigenvalues
$\lm_{\bb}$ of the self-adjoint operator $\bar{A}(\bb,e_0)$ on $\ol{D}(\om,2\pi)$ such that for small
enough $|\beta-\beta_0|$, we have
\be  \bar{A}(\bb,e_0)x_\beta=\lambda_\beta x_\beta,   \lb{4.17}\ee
where $\lambda_{\beta_0}=0$. Taking inner product with $x_\beta$ on both sides of (\ref{4.17})
and then differentiating it with respect to $\beta$ at $\beta_0$, we get
\bea  \frac{\pt}{\pt\bb}\lambda_{\beta}|_{\bb=\bb_0}
&=& \<\frac{\pt}{\pt\bb}\bar{A}(\bb,e_0)x_{\bb},x_{\bb}\>|_{\bb=\bb_0}
    + 2\<\bar{A}(\bb,e_0)x_{\bb},\frac{\pt}{\pt\bb}x_{\bb}\>|_{\bb=\bb_0}  \nn\\
&=& \<\frac{\pt}{\pt\bb}\bar{A}(\bb_0,e_0)x_0,x_0\>  \nn\\
&=& -\frac{1}{(\bb_0+1)^2}\<A(-1,e_0)x_0,x_0\>   \nn\\
&\le& 0, \nn\eea
where the second equality follows from (\ref{4.16}), the last equality follows from the definition
of $\bar{A}(\bb,e)$ and (\ref{4.11}), the last inequality follows from the non-negative definiteness of
$A(-1,e)$ given by Lemma 4.1.
Moreover, assume the last equality holds, then by Lemma 4.1, we must have $\omega=1$ and
\be
x_0=(c_1(1+e\cos t),c_2(1+e\cos t))^T \label{x0}
\ee
for some constant $c_1,c_2\in\C$.
By (\ref{bar.A}), (\ref{4.16}) and (\ref{x0}), we have
\bea
0&=&\<(\frac{A(-1,e)}{\bb_0+1}+\frac{I_2+3S(t)}{2(1+e\cos t)})x_0,x_0\>
\nonumber\\
&=&\<\frac{I_2+3S(t)}{2(1+e\cos t)}x_0,x_0\>
\nonumber\\
&=&\pi(|c_1|^2+|c_2|^2)
\nonumber\\
&>&0,
\eea
where the last inequality follows by $x_0\ne0$.
This is a contradiction.
Thus (\ref{4.15}) is proved. \hb

Consequently we arrive at
\begin{corollary}\label{C4.5} For every fixed $e\in [0,1)$ and $\om\in \U$, the index function
$\phi_{\om}(A(\bb,e))$, and consequently $i_{\om}(\ga_{\bb,e})$, is non-decreasing
as $\bb$ increases from $0$ to $+\infty$.
When $\om=1$, these index functions are increasing and tends from $0$ to $\infty$,
and when $\om\in\U\bs\{1\}$, they are increasing and tends from $2$ to $\infty$.
\end{corollary}

{\bf Proof.} For $0\le \bb_1<\bb_2$ and fixed $e\in [0,1)$, when $\bb$ increases from $\bb_1$ to
$\bb_2$, it is possible that positive eigenvalues of $\bar{A}(\bb_1,e)$ pass through $0$ and become
negative ones of $\bar{A}(\bb_2,e)$, but it is impossible that negative eigenvalues of
$\bar{A}(\bb_2,e)$ pass through $0$ and become positive by (ii) of Lemma \ref{Lemma:increasing.of.index}.
Therefore the first claim holds.

To  prove the second claim, we firstly define a space
\begin{equation}\label{En}
E_n=\span\left\{
R(t)\left(\matrix{0\cr \cos it}\right)\Big|0\le t\le 2\pi,\;i=1,2,...n
\right\}.
\end{equation}
Thus we have $\dim E_n=n$.
Let $\eta(t)$ be a nonzero $C^{\infty}$ function such that $\eta^{(m)}(0)=\eta^{(m)}(2\pi)=0$ for any integer $m\ge0$.
Then we have $\eta(t)E_n\subseteq D(\om,2\pi)$ for any $\om\in\U$.

For any $(\bb,e)\in[0,\infty)\times[0,1)$, $0\ne y(t)=R(t)\left(\matrix{0\cr x(t)}\right)\in E_n$, we have
\begin{eqnarray}
\<A(\bb,e)\eta(t)y(t),\eta(t)y(t)\>
&=&\left\<[-\frac{d^2}{dt^2}I_2-I_2+R(t)K_{\bb,e}(t)R(t)^T]R(t)\left(\matrix{0\cr \eta(t)x(t)}\right),R(t)
\left(\matrix{0\cr \eta(t)x(t)}\right)\right\>
\nonumber\\
&=&\bigg\<[-\frac{d^2R(t)}{dt^2}\left(\matrix{0\cr \eta(t)x(t)}\right)-2\frac{dR(t)}{dt}
\left(\matrix{0\cr (\eta(t)x(t))'}\right)-R(t)\left(\matrix{0\cr (\eta(t)x(t))''}\right)
\nonumber\\
&&+R(t)(-I_2+K_{\bb,e}(t))\left(\matrix{0\cr \eta(t)x(t)}\right)],R(t)\left(\matrix{0\cr \eta(t)x(t)}\right)\bigg\>
\nonumber\\
&=&\bigg\<[R(t)\left(\matrix{0\cr \eta(t)x(t)}\right)-2R(t)J_2\left(\matrix{0\cr (\eta(t)x(t))'}\right)
\nonumber\\
&&+R(t)\left(\matrix{0\cr -(\eta(t)x(t))''-({\bb\over1+e\cos t}+1)(\eta(t)x(t))}\right)],R(t)\left(\matrix{0\cr \eta(t)x(t)}\right)\bigg\>
\nonumber\\
&=&\left\<R(t)\left(\matrix{2(\eta(t)x(t))'\cr -(\eta(t)x(t))''-{\bb\over1+e\cos t}\eta(t)x(t)}\right),R(t)\left(\matrix{0\cr \eta(t)x(t)}\right)\right\>
\nonumber\\
&=&\int_0^{2\pi}[(\eta(t)x(t))']^2dt-\bb\int_0^{2\pi}{(\eta(t)x(t))^2\over1+e\cos t}dt
\nonumber\\
&\le&(C_n-{\bb\over1+e})\int_0^{2\pi}(\eta(t)x(t))^2dt,
\end{eqnarray}
where we have used the property $\eta(t)x(t)|_{t=0}=0$,
and $C_n$ is a constant which depend on space $E_n$ because of the finite dimension of $E_n$.
When $\bb>2C_n>(1+e)C_n$, we obtain that $\<A(\bb,e)\cdot,\cdot\>$ is negative definite on a subspace $\eta(t)E_n$ of $\bar{D}(\om,2\pi)$. Hence
\begin{equation}
i_{\om}(\ga_{\bb,e})\ge n,\quad{\rm if}\;(\bb,e)\in(2C_n,\infty)\times[0,1),
\end{equation}
and together with (\ref{index.of.0e}) on the initial values of index at $\bb=0$, the second part is proved.
\hb

From now on in this section, we will focus on the case of $\omega=1$ and $\omega=-1$. Furthermore, we have

\begin{proposition}\label{P4.6}
When $\bb<\frac{2}{3\sqrt{2}-1}(n^2-\frac{e}{1+e})(1-e)-1$, we have
\begin{equation}
i_1(\ga_{\bb,e})\le4n+2.
\end{equation}
\end{proposition}

{\bf Proof.} Recalling  (\ref{f0}) and (\ref{fn}), for $n\in\N$, we define
\begin{eqnarray}
X_n&=&\span\left\{\pmatrix{1\cr0},\pmatrix{0\cr1}\right\}
\oplus\span\left\{\pmatrix{\cos it\cr0},\pmatrix{0\cr\cos it},
                  \pmatrix{\sin it\cr0},\pmatrix{0\cr\sin it}\Big|i=1,2,...n\right\},
\\
Y_n&=&\span\left\{\pmatrix{\cos it\cr0},\pmatrix{0\cr\cos it},
                  \pmatrix{\sin it\cr0},\pmatrix{0\cr\sin it}\Big|i>n\right\}.
\end{eqnarray}
Then $\bar{D}(1,2\pi)=X_n\oplus Y_n$, $\dim X_n=4n+2$ and $(-\frac{d^2}{dt^2}I_2-I_2)|_{Y_n}\ge n^2-1$.
Moreover, for $y(t)=\pmatrix{y_1(t)\cr y_2(t)}\in  Y_n$, we have
\begin{eqnarray}
&&\int_0^{2\pi}\frac{1}{2(1+e\cos t)}[(\bb+3)I_2+3(\bb+1)S(t)]y(t)\cdot y(t)dt
\nonumber\\
&&=\int_0^{2\pi}\frac{1}{1+e\cos t}y(t)\cdot y(t)dt
+(\bb+1)\int_0^{2\pi}\frac{1}{2(1+e\cos t)}[I_2+3S(t)]y(t)\cdot y(t)dt
\nonumber\\
&&\ge\frac{1}{1+e}|y(t)|_2^2
+(\bb+1)\int_0^{2\pi}\frac{(1+3\cos2t)y_1(t)^2+6\sin2ty_1(t)y_2(t)+(1-3\cos2t)y_2(t)^2}{2(1+e\cos t)}dt
\nonumber\\
&&\ge\frac{1}{1+e}|y(t)|_2^2
+(\bb+1)\int_0^{2\pi}\frac{(1+3\cos2t-3|\sin2t|)y_1(t)^2+(1-3\cos2t-3|\sin2t|)y_2(t)^2}{2(1+e\cos t)}dt
\nonumber\\
&&\ge\frac{1}{1+e}|y(t)|_2^2
+(\bb+1)\int_0^{2\pi}\frac{(1-3\sqrt{2})y_1(t)^2+(1-3\sqrt{2})y_2(t)^2}{2(1+e\cos t)}dt
\nonumber\\
&&\ge(\frac{1}{1+e}-\frac{(3\sqrt{2}-1)(\bb+1)}{2(1-e)}||y||_2^2
\end{eqnarray}
Thus for any $y(t)=\pmatrix{y_1(t)\cr y_2(t)}\in  Y_n$, we obtain
\begin{equation}
\<A(\bb,e)y(t),y(t)\>\ge(n^2-1+\frac{1}{1+e}-\frac{(3\sqrt{2}-1)(\bb+1)}{2(1-e)})||y||_2^2,
\end{equation}
and hence when $\bb<\frac{2}{3\sqrt{2}-1}(n^2-\frac{e}{1+e})(1-e)-1$, we have $\<A(\bb,e)y(t),y(t)\>\ge0$
for any $y(t)\in Y_n$.
Then it implies $i_1(\ga_{\bb,e})\le4n+2$.
\hb


\subsection{The degenerate curves of elliptic Euler solution}

Because $A(\bb,e)$ is a self-adjoint operator on $\bar{D}(\om,2\pi)$,
and a bounded perturbation of the operator $-\frac{d^2}{dt^2}I_2$, then $A(\bb,e)$
has discrete spectrum on $\bar{D}(\om,2\pi)$.
Thus we can define the $n$-th degenerate point for any $\om$ and $e$:
\be\lb{bb_s}
\bb_n(\om,e)=\min\left\{\bb>0\;\bigg|\;
  \begin{array}{l}
   [i_\om(\ga_{\bb,e})+v_\om(\ga_{\bb,e})]-[i_\om(\ga_{0,e})+v_\om(\ga_{0,e})]\ge n
    \end{array}\right\}.
\ee
By Lemma \ref{Lemma:increasing.of.index} $(iii)$,
$i_\om(\ga_{\bb,e})+v_\om(\ga_{\bb,e})$ is a right continuous step function with respect to $\bb$.
Additionally, by Corollary \ref{C4.5}, $i_\om(\ga_{\bb,e})+v_\om(\ga_{\bb,e})$ tends to $+\infty$ as
$\bb\rightarrow+\infty$, the minimum of the right hand side in (\ref{bb_s}) can be obtained.
Indeed, $\ga_{\bb,e}$ is $\om$-degenerate at point $(\bb_n(\om,e),e)$, i.e.,
\be\label{degenerate.of.bn}
\nu_\om(\ga_{\bb_n(\om,e),e})\ge 1.
\ee
Otherwise, if there existed some small enough $\ep>0$ such that $\bb=\bb_n(\om,e)-\ep$ would
satisfy $[i_\om(\ga_{\bb,e})+v_\om(\ga_{\bb,e})]-[i_\om(\ga_{0,e})+v_\om(\ga_{0,e})]\ge n$ in (\ref{bb_s}),
it would yield a contradiction.

For fixed $\om$ and $n$, $\bb_n(\om,e)$ actually forms a curve with respect to the eccentricity $e\in [0,1)$
as we shall prove below in this section, which we called the $n$-th $\om$-degenerate curve. By Corollary
\ref{C4.5}, $\bb_n(\om,e)$ is non-decreasing with respect to $n$ for fixed $\om$ and $e$.
We have

\begin{lemma}\lb{Lemma4.5}
For any fixed $n\in\N$ and $\om\in\U$, the degenerate curve $\bb_n(\om,e)$ is continuous with respect to $e\in[0,1)$.
\end{lemma}

{\bf Proof.} In fact, if the function $\bb_n(\om,e)$ is not continuous in $e\in[0,1)$,
then there exists some $\tilde{e}\in[0,1)$, a sequence $\{e_i|i\in\N\}\subset[0,1)\backslash\{\tilde{e}\}$
and $\bb_0\ge0$ such that
\be\label{converge.to.bb0}
\bb_n(\om,e_i)\to\bb_0\ne\bb_n(\om,\tilde{e})   \quad{\rm and}\quad e_i\to\tilde{e}  \quad{\rm as}\quad i\to+\infty.
\ee
By (\ref{degenerate.of.bn}), we have $\om\in\sigma(\ga_{\bb_n(\om,e_i),e_i}(2\pi))$.
By the continuity of eigenvalues of $\ga_{\bb_n(\om,e_i),e_i}(2\pi)$ in $e_i$ as $i\to +\infty$ and
(\ref{converge.to.bb0}), we have $\om\in\sigma(\ga_{\bb_0,\tilde{e}}(2\pi))$,
and hence
\be\label{bb0.om-degenerate}
\nu_\om(\ga_{\bb_0,\tilde{e}})\ge1.
\ee
We continue in two cases according to the sign of the difference $\bb_0-\bb_n(\om,\tilde{e})$.
For convenience, let
\be
g(\bb,e)=[i_\om(\ga_{\bb,e})+v_\om(\ga_{\bb,e})]-[i_\om(\ga_{0,e})+v_\om(\ga_{0,e})].
\ee

If $\bb_0<\bb_n(\om,\tilde{e})$, firstly we must have $g(\bb_0,\tilde{e})<n$,
otherwise by the definition of $\bb_n(\om,\tilde{e})$,  we must have $\bb_n(\om,\tilde{e})\le\bb_0$.

Let $\tilde\bb\in(\bb_0,\bb_n(\om,\tilde{e}))$ such that $\nu_\om(\ga_{\bb,\tilde{e}})=0$
for any $\bb\in(\bb_0,\tilde\bb]$. By the continuity of eigenvalues of $\ga_{\bb,e}(2\pi)$ with
respect to $\bb$ and $e$, there exists a neighborhood $\mathcal{O}$ of $(\tilde\bb,\tilde{e})$
such that $\nu(\ga_{\bb,e})=0$ for any $(\bb,e)\in\mathcal{O}$. Then $i_\om(\ga_{\bb,e})$, and
hence $g(\bb,e)$ is constant in $\mathcal{O}$. By (\ref{converge.to.bb0}), for $i$ large enough,
we have $\bb_n(\om,e_i)<\tilde\bb$ and $(\tilde\bb,e_i)\in\mathcal{O}$, and hence
$g(\tilde\bb,e_i)\ge g(\bb_n(\om,e_i),e_i)\ge n$. Therefore, we have $g(\tilde\bb,\tilde{e})\ge n$.
By the definition of (\ref{bb_s}), we have $\bb_n(\om,\tilde{e})\le\tilde\bb$
which contradicts $\tilde\bb\in(\bb_0,\bb_n(\om,\tilde{e}))$.

If $\bb_0>\bb_n(\om,\tilde{e})$,
there exists $\bar\bb\in(\bb_n(\om,\tilde{e}),\bb_0)$ such that $\nu_\om(\ga_{\bb,\tilde{e}})=0$
for any $\bb\in(\bb_n(\om,\tilde{e}),\bar\bb]$.
By the continuity of eigenvalues of $\ga_{\bb,e}(2\pi)$ with respect to $\bb,e$,
there exists a neighborhood $\mathcal{U}$ of $(\bar\bb,\tilde{e})$ such that $\nu(\ga_{\bb,e})=0$
for any $(\bb,e)\in\mathcal{U}$.
Then $i_\om(\ga_{\bb,e})$, and hence $g(\bb,e)$ is constant in $\mathcal{U}$.
By (\ref{converge.to.bb0}), for $i$ large enough,
we have $\bar\bb<\bb_n(\om,e_i)$ and $(\bar\bb,e_i)\in\mathcal{U}$.
$g(\bar\bb,e_i)=g(\bar\bb,\tilde{e})\ge n$ implies $\bb_n(\om,e_i)\le \bar\bb$, a contradiction.

Thus the continuity of $\bb_n(\om,e)$ in $e\in[0,1)$ is proved.
\hb

For $n=1$, by
Corollary \ref{C4.5}, we have another equivalent definition:
\be\lb{bb_1}
\bb_1(\om,e)=\min\{\bb>0\;|\;A(\bb,e)\; is\;degenerate\;on\;\bar{D}(\om,2\pi)\}.
\ee
Moreover, let $\om=1$, we have the following theorem

\begin{theorem}\label{Th:near.bb=0}
For any $\epsilon>0$, there exists a $\bb_0=\bb_0(\epsilon)>0$ such that
\be
\bb_1(1,e)>\bb_0,\quad\forall e\in[0,1-\epsilon]. \lb{bb_s.g.0}
\ee
\end{theorem}

{\bf Proof.} By the fact that $A(\bb,e)$ has discrete spectrum and definition (\ref{bb_s}) we have
$\bb_1(1,e)>0$ for fixed $e\in[0,1)$.
If (\ref{bb_s.g.0}) does not hold, there is an sequence $\{e_n\}_{n=1}^{\infty}\subseteq[0,1-\epsilon]$
such that $\lim_{n\to\infty}e_n=e_0$ for some $e_0\in[0,1-\epsilon]$ and $\lim_{n\to\infty}\bb_1(1,e_n)=0$.
We consider the operator $A(\frac{1}{2}\bb_1(1,e_0),e_0)$. It is non-degenerate by the definition of
$\bb_1(1,e)$ in (\ref{bb_1}). Therefore, $A(\bb,e)$ is non-degenerate and has the same indices with
$A(\frac{1}{2}\bb_1(1,e_0),e_0)$, when $(\bb,e)$ is in a small neighborhood of
$(\frac{1}{2}\bb_1(1,e_0),e_0))$. Moreover, $\phi_1(A(\frac{1}{2}\bb_1(1,e_0),e_0))=v(A(0,e_0))=3$ by
Lemma 4.2. Then for $n$ large enough we obtain
$$  \phi_1(A(\frac{1}{2}\bb_1(1,e_0),e_n)) \;=\; \phi_1(A(\frac{1}{2}\bb_1(1,e_0),e_0)) \;=\; 3.   $$

On the other hand, by the non-decreasing property of $i_1(A(\bb,e))$ with respect to $\bb$, and notice that
$v(A(\bb_1(1,e_n),e_n))\ge1$ by definition (\ref{bb_s}),
for $n$ sufficiently large, we have ${1\over2}\bb_1(1,e_0)>\bb_1(1,e_n)$ and
\bea \phi_1(A(\frac{1}{2}\bb_1(1,e_0),e_n))
&\ge& \phi_1(A(\bb_1(1,e_n),e_n)) + v_1(A(\bb_1(1,e_n),e_n)) \nn\\
&\ge& \phi_1(A(0,e_n))+1 \nn\\
&\ge& 4. \eea
where we have applied (\ref{2.30}) and Lemma 4.2 (iii). This is a contradiction.
Thus the theorem is proved.
\hb

We now calculate the intersection points of the $1$-degenerate curves with the horizontal axis.
Recall (\ref{A2}) and (\ref{A3}), for $\hat{\bb}_n$ defined by (\ref{om_n}), $A(\hat{\bb}_n,0)$
is degenerate and
\begin{equation}
\ker A(\hat{\bb}_n,0)=\span\left\{
R(t)\left(\matrix{a_n\sin nt\cr \cos nt}\right),\quad
R(t)\left(\matrix{a_n\cos nt\cr -\sin nt}\right)
\right\},
\end{equation}
where $a_n=\frac{n^2-\hat\bb_n}{2n}$.

\begin{remark}
By (\ref{A3}), $A(\bb,0)R(t)\left(\matrix{a_n\sin nt\cr \cos nt}\right)=0$ reads
\be
\left\{
\begin{array}{cr}
n^2a_n-2n+(2\bb+3)a_n&=0,
\\
n^2-2na_n-\bb&=0.
\end{array}
\right.
\ee
Then $2\bb^2-(n^2-3)\bb-n^2(n^2-1)=0$ which yields $\bb=\hat\bb_n$ again and $a_n=\frac{n^2-\hat\bb_n}{2n}$.
Moreover, by (\ref{bb_n.approx}), $a_n\approx\frac{n^2-(n^2-4/3)}{2n}=\frac{2}{3n}$.
\end{remark}

Thus every $1$-degenerate curve starts from the point $(\hat\bb_n,0)$.
Moreover we have
\begin{lemma}\label{L4.9}
\begin{equation}
\bb_n(1,0) = \hat\bb_{m+1}, \quad {\rm if}\;\;n = 2m-1\;\;{\rm or}\;\;2m. \lb{bb_n0}
\end{equation}
\end{lemma}

{\bf Proof.} By (\ref{1-index.of.b0}) and (\ref{null.1-index.of.b0}), we have $i_1(0,0)+v_1(0,0)=3$,
$v_1(\hat\bb_{m+1},0)=2$ and
\begin{equation}
[i_1(\bb,0)+v_1(\bb,0)]-[i_1(0,0)+v_1(0,0)]
\left\{
\begin{array}{l}
\le 2m-2, \quad {\rm if}\;\;\bb<\hat\bb_{m+1},
\\
=2m, \quad {\rm if}\;\;\hat\bb_{m+1}\le\bb<\hat\bb_{m+2},
\\
\ge 2m+2, \quad {\rm if}\;\;\bb\ge\hat\bb_{m+2}.
\end{array}
\right.
\end{equation}
For $n=2m-1$ or $2m$, $[i_1(\bb,0)+v_1(\bb,0)]-[i_1(0,0)+v_1(0,0)]\ge n$ is equivalent to $\bb\ge\hat\bb_{m+1}$.
Then the minimal value of $\bb$ in $\{\bb\ge\hat{\bb}_{m+1}\}$ such that $A(\bb,e)$ is degenerate on
$\overline{D}(1,2\pi)$ is $\hat{\bb}_{m+1}$. Thus by (\ref{bb_s}), we obtain (\ref{bb_n0}). \hb

Moreover, we have the following theorem:
\begin{theorem}\label{Th:multiplicity}
Every $1$-degenerate curves has even multiplicity.
\end{theorem}

{\bf Proof.} The statement has already been proved for $e=0$. We will prove that, if $A(\bb,e)z=0$
has a solution $z\in\bar{D}(1,2\pi)$ for a fixed value $e\in (0,1)$, there exists a second periodic
solution which is independent of $z$. Then the space of solutions of $A(\bb,e)z=0$ is the direct
sum of two isomorphic subspaces, hence it has even dimension. This method is due to R. Mat\'{i}nez,
A. Sam\`{a} and C. Sim\`{o} in \cite{MSS1}.

Let $z(t)=R(t)(x(t),y(t))^T$ be a nontrivial solution of $A(\bb,e)z(t)=0$, then it yields
\begin{equation}
\left\{
\begin{array}{l}
(1+e\cos t)x''(t)=(2\bb+3)x(t)+2y'(t)(1+e\cos t),
\\
(1+e\cos t)y''(t)=-\bb y(t)-2x'(t)(1+e\cos t).
\end{array}
\right.
\end{equation}
By Fourier expansion, $x(t)$ and $y(t)$ can be written as
\bea
x(t) =a_0+\sum_{n\ge1}a_n\cos nt+\sum_{n\ge1}b_n\sin nt,
\\
y(t) =c_0+\sum_{n\ge1}c_n\cos nt+\sum_{n\ge1}d_n\sin nt.
\eea
Then the coefficient must satisfy the following uncoupled sets of recurrences:
\begin{equation}
\left\{
\begin{array}{l}
(2\bb+3)a_0=-e(d_1+\frac{a_1}{2}),
\\
eA_2\pmatrix{a_2\cr d_2}=B_1\pmatrix{a_1\cr d_1},
\\
eA_{n+1}\pmatrix{a_{n+1}\cr d_{n+1}}
=B_n\pmatrix{a_n\cr d_n}-eA_{n-1}\pmatrix{a_{n-1}\cr d_{n-1}},
\quad n\ge2,
\end{array}
\right.
\label{ad.equations}
\end{equation}
and
\begin{equation}
\label{bc.equations}
\left\{
\begin{array}{l}
-\bb c_0=e(b_1-\frac{c_1}{2}),
\\
eA_2\pmatrix{b_2\cr -c_2}=B_1\pmatrix{b_1\cr -c_1},
\\
eA_{n+1}\pmatrix{b_{n+1}\cr -c_{n+1}}
=B_n\pmatrix{b_n\cr -c_n}-eA_{n-1}\pmatrix{b_{n-1}\cr -c_{n-1}},
\quad n\ge2,
\end{array}
\right.
\end{equation}
where
\begin{equation}
A_n=-\frac{n}{2}\pmatrix{n&2\cr 2&n},\quad
B_n=\pmatrix{n^2+2\bb+3&2n\cr 2n&n^2-\bb}.
\end{equation}

Thus $\det(B_1)=-2\bb(\bb+1)\ne0$ for $\bb>0$ and $\det(A_n)\ne0$ when $n\ne3$.
Thus given $(a_2,d_2)^T$, we can obtain $(a_1,d_1)^T$ uniquely from the second equality of (\ref{ad.equations}),
and then obtain $(a_n,d_n)^T$ for $n\ge3$ by the last equality of (\ref{ad.equations}).

By the non-triviality of $z=z(t)$, both (\ref{ad.equations}) and (\ref{bc.equations}) have solutions
$\{(a_n,d_n)\}_{n=1}^\infty$ and $\{(b_n,c_n)\}_{n=1}^\infty$ respectively.
We assume (\ref{ad.equations}) admits a nontrivial solutions.
Then $\sum_{n\ge1}a_n\cos nt$ and $\sum_{n\ge1}d_n\sin nt$ are convergent.
Thus, $\sum_{n\ge1}a_n\sin nt$ and $-\sum_{n\ge1}d_n\cos nt$ are convergent too.
Moreover, by the similar structure between equations (\ref{ad.equations}) and (\ref{bc.equations}),
we can construct a new solution of ((\ref{bc.equations})) given below
\bea
\tilde{c}_0&=&-\frac{e}{\bb}(a_1+\frac{d_1}{2}),
\\
\left(\matrix{\tilde{b}_n\cr\tilde{c}_n}\right)&=&\left(\matrix{a_n\cr -d_n}\right),\quad n\ge1.
\eea
Therefore we can build two independent solutions of $A(\bb,e)w=0$ as
\bea
w_1(t)&=&R(t)\pmatrix{a_0+\sum_{n\ge1}a_n\cos nt\cr\sum_{n\ge1}d_n\sin nt},
\\
w_2(t)&=&R(t)\pmatrix{\sum_{n\ge1}\tilde{b}_n\sin nt\cr \tilde{b}_0+\sum_{n\ge1}\tilde{c}_n\cos nt}
        =R(t)\pmatrix{\sum_{n\ge1}a_n\sin nt\cr -\frac{e}{\bb}(a_1+\frac{d_1}{2})-\sum_{n\ge1}d_n\cos nt}.
\eea
\hb

\begin{remark}
In the above proof, if $b_n=\lambda a_n$, $c_n=-\lambda d_n$ for $n\ge1$ and some $\lambda\ne0$,
we can construct two independent solutions.
But if this situation does not hold,
and both $(a_n, d_n)^T$, $(b_n, c_n)^T$ are nontrivial sequences,
then we can construct four independent solutions by the similar method.
In the following Remark \ref{Remark:multiplicity.2}, we will show that the latter situation does not appear.
\end{remark}

\begin{theorem}\label{Th:odd.indices}
For any $\bb>0$ and $0\le e<1$, $i_1(\ga_{\bb,e})$ is an odd number.
\end{theorem}
{\bf Proof.} When $e=0$, the conclusion of our theorem follows from (\ref{1-index.of.b0}).

Now we suppose $0<e<1$. By Lemma \ref{Lemma:increasing.of.index} $(iii)$,
we can choose an $\epsilon_0>0$ small enough such that for any $\epsilon\in (0,\epsilon_0)$, by
(\ref{index.of.0e}) and (\ref{null.index.of.0e}) we obtain
\begin{equation}
i_1(\ga_{\epsilon,0})=i_1(\ga_{0,e})+\nu_1(\ga_{0,e})=3.
\end{equation}

Now for any $\bb_{\ast}\ge \frac{\epsilon_0}{2}$,
by Lemma \ref{Lemma:increasing.of.index} and Corollary \ref{C4.5},
the set
$\{\frac{\epsilon_0}{2}<\bb\le\bb_{\ast}\;|\;\nu_1(\ga_{\bb,e})\ne 0\}$ contains only finitely many points.
Thus we can suppose
\begin{equation}
\{\frac{\epsilon_0}{2}\le\bb\le\bb_{\ast}\;|\;\nu_1(\ga_{\bb,e})\ne 0\}=\{\bb_{\ast 1},\ldots,\bb_{\ast n}\}.
\end{equation}
Then by Lemma \ref{Lemma:increasing.of.index} $(iii)$, we have
\begin{equation}
\label{odd.sum}
i_1(\ga_{\bb_{\ast},0}) = i_1(\ga_{\ep_0/2,e})+\sum_{k=1}^n\nu_1(\ga_{\bb_{\ast k},e})
    = 3+\sum_{k=1}^n\nu_1(\ga_{\bb_{\ast k},e}).
\end{equation}
By the proof of Theorem \ref{Th:multiplicity} and its remark, every $\nu_1(\ga_{\bb_{\ast k},e})$ is even for $1\le k\le n$.
Thus $i_1(\ga_{\bb_{\ast},0})$ is odd by (\ref{odd.sum}).
\hb

\subsection{The order of the degenerate curves and the normal forms of $\ga_{\bb,e}(2\pi)$}

Now we study the order of the $1$-degenerate curves and $-1$-degenerate curves.

\begin{theorem}\label{Th:no.intersection}
Any $1$-degenerate curve and any $-1$-degenerate curve cannot intersect each other.
That is, for any $0<e<1$, there does not exist $n_1,n_2\in\N$ such that $\bb_{n_1}(1,e)=\bb_{n_2}(-1,e)$.
\end{theorem}

{\bf Proof. } If not, suppose $(\bb_{\ast},e_{\ast})$ with $\bb_{\ast}>0$ and $0<e_{\ast}<1$ is an
intersection point of some $1$-degenerate curve and a $-1$-degenerate curve.
Then $\nu_1(\ga_{\bb_{\ast},e_{\ast}})\ge1$ and $\nu_{-1}(\ga_{\bb_{\ast},e_{\ast}})\ge1$.
Moreover, by Theorem \ref{Th:orth} and its remark, $\nu_1(\ga_{\bb_{\ast},e_{\ast}})\ge1$ is even.
Therefore, there exists a $b\in\R$ such that $\ga_{\bb_{\ast},e_{\ast}}(2\pi)\in \Sp(4)$ satisfies:
\begin{equation}
\ga_{\bb_{\ast},e_{\ast}}(2\pi)\approx I_2\diamond N_1(-1,b).
\end{equation}

By Lemma \ref{Lm:Path.sum}, there exist two paths $\ga_i\in\P_{2\pi}(2)$ such that we have
$\ga_1(2\pi)=I_2$, $\ga_2(2\pi)=N_1(-1,b)$, $\ga_{\bb_{\ast},e_{\ast}}\sim\ga_1\dm\ga_2$, and
$i_1(\ga_{\bb_{\ast},e_{\ast}})=i_1(\ga_1)+i_1(\ga_2)$. By Theorem 8.1.4 and Theorem 8.1.5 on
pp.179-181 of \cite{Lon4}, both $i_1(\ga_1)$ and $i_1(\ga_2)$ must be odd numbers. Therefore
$i_1(\ga_{\bb_{\ast},e_{\ast}})$ must be even. But Theorem \ref{Th:odd.indices} tell us
$i_1(\ga_{\bb_{\ast},e_{\ast}})$ is an odd number. It is a contradiction.
\hb

Because of the starting points from $\bb$-axis of the $1$-degenerate curves and $-1$-degenerate curves
are alternatively distributed, and these curves are analytic by Theorem \ref{Th:analytic} and Theorem
\ref{Th:analytic.om}, any two $1$-degenerate curves (or two $-1$-degenerate curves) starting from
different points cannot intersect each other. Thus we have the following corollary:
\begin{corollary}\label{Cor:order}
Using notations in Theorem \ref{T1.3}, the $1$-degenerate curves and $-1$-degenerate curves of the
elliptic Euler solutions in Figure 1 can be ordered from left to right by
\begin{equation}
0,\; \Xi_1^-,\; \Xi_1^+,\; \Ga_1,\; \Xi_2^-,\; \Xi_2^+,\; \Ga_2,\; \ldots,\; \Xi_n^-,\; \Xi_n^+,\; \Ga_n.
\end{equation}
More precisely, for each fixed $e\in [0,1)$, we have
\begin{eqnarray}
0&<&\bb_1(-1,e)\le \bb_2(-1,e)<\bb_1(1,e)=\bb_2(1,e)<\bb_3(-1,e)\le \bb_4(-1,e)<\bb_3(1,e)=\bb_4(1,e)<\cdots
\nonumber
\\
&<&\bb_{2m-1}(-1,e)\le \bb_{2m}(-1,e)<\bb_{2m-1}(1,e)=\bb_{2m}(1,e)<\cdots
\end{eqnarray}
\end{corollary}

\begin{remark}\label{Remark:multiplicity.2}
By Theorem \ref{Th:orth}, Theorem \ref{Th:multiplicity} and (\ref{null.1-index.of.b0}),
the $1$-degenerate curves start form $(\hat\bb_n,0)$ with multiplicity $2$ near $e=0$.
If there is some point $(\bb_0,e_0)\in(0,+\infty)\times(0,1)$ such that $\nu_{1}(\ga_{\bb_0,e_0})\ge4$.
Then there must exist two different $1$-degenerate curves which intersect at $(\bb_0,e_0)$.
This contradicts Corollary \ref{Cor:order}.
Thus every $1$-degenerate curve has exact multiplicity $2$.
\end{remark}

By a similar proof of Theorem \ref{Th:no.intersection}, we have
\begin{theorem}\label{Th:no.intersection2}
For $\om\ne\pm1$, any $\om$-degenerate curve and any $-1$-degenerate curve cannot intersect each other.
That is, for any $0<e<1$, there does not exist $n_1,n_2\in\N$ such that $\bb_{n_1}(\om,e)=\bb_{n_2}(-1,e)$.
\end{theorem}

Now we can give

{\bf The Proof of Theorem 1.5.} (i) follows from the discussion on (46) of \cite{HS2}.

(ii) If $0<\bb<\bb_1(-1,e)$, then by the definitions of the degenerate curves and
Lemma \ref{Lemma:increasing.of.index} $(iii)$, we have
\begin{equation}
i_1(\ga_{\bb,e})=3,\quad\quad \nu_1(\ga_{\bb,e})=0, \lb{4.89}
\end{equation}
and
\begin{equation}
i_{-1}(\ga_{\bb,e})=2,\quad\quad \nu_{-1}(\ga_{\bb,e})=0.  \lb{4.90}
\end{equation}

Firstly, if $\ga_{\bb,e}(2\pi)\approx N_2(e^{\sqrt{-1}\theta},b)$ for some $\theta\in(0,\pi)\cup(\pi,2\pi)$,
we have
\be
i_{-1}(\ga_{\bb,e})=i_1(\ga_{\bb,e})-S_{N_2(e^{\sqrt{-1}\theta},b)}^-(e^{\sqrt{-1}\theta})
     +S_{N_2(e^{\sqrt{-1}\theta},b)}^+(e^{\sqrt{-1}\theta})=i_1(\ga_{\bb,e})
\ee
or
\be
i_{-1}(\ga_{\bb,e})=i_1(\ga_{\bb,e})-S_{N_2(e^{\sqrt{-1}\theta},b)}^-(e^{\sqrt{-1}(2\pi-\theta)})
    +S_{N_2(e^{\sqrt{-1}\theta},b)}^+(e^{\sqrt{-1}(2\pi-\theta)})=i_1(\ga_{\bb,e}),
\ee
which contradicts to (\ref{4.89}) and (\ref{4.90}).

Then we can suppose $\ga_{\bb,e}(2\pi)\approx M_1\diamond M_2$ where $M_1$ and $M_2$ are two basic normal
forms in $\Sp(2)$ defined in Section 5.2 below. By Lemma \ref{Lm:Path.sum} there exist two paths $\ga_1$ and
$\ga_2$ in $\P_{2\pi}(2)$ such that $\ga_1(2\pi)=M_1$, $\ga_2(2\pi)=M_2$, $\ga_{\bb,e}\sim\ga_1\dm\ga_2$,
and $i_1(\ga_{\bb,e})=i_1(\ga_1)+i_1(\ga_2)$ hold.

Thus one of $i_1(\ga_1)$ and $i_2(\ga_1)$ must be odd, and the other is even. Without loss of generality, we
suppose $i_1(\ga_2)$ is odd. Notice that $\nu_1(\ga_{\bb,e})=0$, by Theorems 8.1.4 to 8.1.7 on pp.179-183 of
\cite{Lon4} and using notations there, we must have $M_2\in\Sp^{th}(2)$ and $\alpha(M_2)=0$. Therefore, $M_2=D(2)$.
Using the same method, we have $M_1=D(-2)$ or $M_1=R(\theta)$ for some $\theta\in(0,\pi)\cup(\pi,2\pi)$.
If $M_1=D(-2)$, by the properties of splitting numbers in Chapter 9 of \cite{Lon4}, specially (9.3.3) on p.204, we
obtain $i_{-1}(\ga_{\bb,e})=i_1(\ga_{\bb,e})$, which contradicts to (\ref{4.89}) and (\ref{4.90}). Therefore, we must
have $M_1=R(\theta)$.

If $\theta\in(0,\pi)$, we have
$i_{-1}(\ga_{\bb,e})=i_1(\ga_{\bb,e})-S_{R(\theta)}^-(e^{\sqrt{-1}\theta})+S_{R(\theta)}^+(e^{\sqrt{-1}\theta})=2$.
When $\theta\in(\pi,2\pi)$, we obtain
$i_{-1}(\ga_{\bb,e})=i_1(\ga_{\bb,e})-S_{R(\theta)}^-(e^{\sqrt{-1}(2\pi-\theta)})+S_{R(\theta)}^+(e^{\sqrt{-1}(2\pi-\theta)})=4$.
Therefore, we have $\theta\in(0,\pi)$, and then $\ga_{\bb,e}(2\pi)\approx R(\theta)\diamond D(2)$. Thus (ii) is proved.

(v) If $\bb_1(-1,e)\ne\bb_2(-1,e)$ and $\bb_1(-1,e)<\bb<\bb_2(-1,e)$,
then by the definitions of the degenerate curves and
Lemma \ref{Lemma:increasing.of.index} $(iii)$, we have
\begin{equation}
i_1(\ga_{\bb,e})=3,\quad\quad \nu_1(\ga_{\bb,e})=0, \lb{4.94}
\end{equation}
and
\begin{equation}
i_{-1}(\ga_{\bb,e})=3,\quad\quad \nu_{-1}(\ga_{\bb,e})=0.  \lb{4.95}
\end{equation}

If $\ga_{\bb,e}(2\pi)\approx N_2(e^{\sqrt{-1}\theta},b)$ in Subsection 5.2 for some
$\theta\in(0,\pi)\cup(\pi,2\pi)$, we now cannot use the method in (ii) directly to obtain the contradiction
because of $i_1(\ga_{\bb,e})=i_{-1}(\ga_{\bb,e})$.

On the one hand, $\ga_{\bb,e}(2\pi)\approx N_2(e^{\sqrt{-1}\theta},b)$ implies that
$(\bb,e)$ is on some $\om$-degenerate curve $\Th_{\om}$ where $\om\ne\pm1$.
On the other hand, $\bb_1(-1,e)<\bb<\bb_2(-1,e)$ implies that
$(\bb,e)$ is between the two $-1$-degenerate curves $\Xi_1^{\pm}$ which
start from the same point $(\hat\bb_{3\over 2},0)$. But $\Th_{\om}$ is a continuous curve defined on
the closed interval $[0,1)$ by Lemma \ref{Lemma4.5}. Thus $\Th_{\om}$ must come down from the point
$(\bb,e)$ to the horizontal axis of $e=0$, and then it must intersect with at least one of $\Xi_1^{\pm}$,
which contradicts Theorem \ref{Th:no.intersection2}.

Then we can suppose $\ga_{\bb,e}(2\pi)\approx M_1\diamond M_2$, and following a similar steps in (ii),
we can obtain $\ga_{\bb,e}(2\pi)\approx D(-2)\diamond D(2)$.

By the same method, (iii)-(iv) and (vi)-(xiv) can be proved and the details is thus omitted here. \hb

\subsection{The two $\omega=1$ degenerate curves coincide and orthogonal to the horizontal axis}

Recall $A(-1,e)$ is non-negative definite on $\overline{D}(1,2\pi)$, and (\ref{Ker.A.-1.e}) holds.
Let $P_1(e)$ be the projection operator from $\overline{D}(1,2\pi)$ to $\ker A(-1,e)$,
then $A(-1,e)+P_1(e)$ is positive definite on its domain $\overline{D}(1,2\pi)$.
Now we set
\be
B(\bb,e)=[A(-1,e)+P_1(e)]^{-\frac{1}{2}}
         \left(\frac{I_2+3S(t)}{2(1+e\cos t)}-\frac{P_1(e)}{\bb+1}\right)
         [A(-1,e)+P_1(e)]^{-\frac{1}{2}}.
\ee
Then we have
\begin{lemma}\label{L4.6}
For $0\le e<1$, $A(\bb,e)$ is $1$-degenerate if and only if $-\frac{1}{\bb+1}$ is an eigenvalue of $B(\bb,e)$.
\end{lemma}
{\bf Proof.} Suppose $A(\bb,e)x=0$ holds for some $x\in\overline{D}(1,2\pi)$.
Let $y=[A(-1,e)+P_1(e)]^{\frac{1}{2}}x$. Then by (\ref{4.11}) we obtain
\bea
&&[A(-1,e)+P_1(e)]^{\frac{1}{2}}\left(\frac{1}{\bb+1}+B(\bb,e)\right)y(t)
\nonumber\\
&&=\left(\frac{A(-1,e)+P_1(e)}{\bb+1}+\frac{I_2+3S(t)}{2(1+e\cos t)}-\frac{P_1(e)}{\bb+1}\right)x(t)
\nonumber\\
&&=\frac{1}{\bb+1}A(\bb,e)x
\nonumber\\
&&=0.    \label{A.B}
\eea

Conversely, if $(\frac{1}{\bb+1}+B(\bb,e))y=0$,
then $x=[A(-1,e)+P_1(e)]^{-\frac{1}{2}}y$ is an eigenfunction of $A(\bb,e)$ belonging to the eigenvalue $0$
by our computations (\ref{A.B}).
\hb

Although $e<0$ does not have physical meaning, we can extend the fundamental solution to the case $e\in(-1,1)$
mathematically and all the above results which holds for $e>0$ also holds for $e<0$.
Then we have
\begin{theorem}\label{Th:analytic}
Every $1$-degenerate curve $(\bb_n(1,e),e)$ in $e\in(-1,1)$ is a real analytic function.
\end{theorem}

{\bf Proof.}
By Lemma \ref{L4.6}, $-\frac{1}{\bb_i(1,e)+1}$ is an eigenvalue of $B(\bb,e)$.
Note that $B(\bb,e)$ is a compact operator and self adjoint when $\bb, e$ are real.
Moreover, it depends analytically on $\bb$ and $e$, and we denote its eigenvalue by $f(\bb,e)$.
By \cite{Ka}(Theorem 3.9 in p.392), we know that $-\frac{1}{\bb_i(1,e)+1}$ is analytical in $e$ for each $i\in\N$.
By Theorem \ref{Th:multiplicity}, Corolary \ref{Cor:order} and Remark \ref{Remark:multiplicity.2},
every $1$-degenerate curve has multiplicity 2,
and any two different $1$-degenerate curves cannot intersect each other.
We can suppose
\be\label{f.be}
-\frac{1}{\bb_i(1,e)+1}=f(\bb_i(1,e),e).
\ee
Differentiate $B(\bb,e)$ with respect to $\bb$, we obtain
\be
\frac{\partial B(\bb,e)}{\partial\bb}=\frac{1}{(\bb+1)^2}[A(-1,e)+P_1(e)]^{-\frac{1}{2}}
         P_1(e)[A(-1,e)+P_1(e)]^{-\frac{1}{2}}
>0.
\ee
By the same techniques in the proof of Lemma \ref{Lemma:increasing.of.index} $(i)$,
we can choose a smooth path of unit norm eigenvectors $x_{\bb,e}$
belongs to a smooth path of real eigenvalues $f(\bb,e)$
of the self adjoint operator $B(\bb,e)$ on $\overline{D}(1,2\pi)$,
it yields
\bea
\frac{\partial f(\bb,e)}{\partial\bb}&=&\<\frac{\partial B(\bb,e)}{\partial\bb}x_{\bb,e},x_{\bb,e}\> \nonumber
\\
&=&\<\frac{1}{(\bb+1)^2}[A(-1,e)+P_1(e)]^{-\frac{1}{2}}P_1(e)[A(-1,e)+P_1(e)]^{-\frac{1}{2}}x_{\bb,e},x_{\bb,e}\>
\nonumber
\\
&\le&\<\frac{1}{(\bb+1)^2}[A(-1,e)+P_1(e)]^{-\frac{1}{2}}(A(-1,e)+P_1(e))[A(-1,e)+P_1(e)]^{-\frac{1}{2}}x_{\bb,e},x_{\bb,e}\>
\nonumber
\\
&=&\frac{1}{(\bb+1)^2}\<x_{\bb,e},x_{\bb,e}\> \nonumber
\\
&=&\frac{1}{(\bb+1)^2},
\eea
where the third equality holds for some $(\bb_0,e_0)\in(0,\infty)\times(-1,1)$
if and only if there exists a nontrivial $x_{\bb_0,e_0}$ such that
\be\label{Bx.x=0}
\<\frac{1}{(\bb+1)^2}[A(-1,e_0)+P_1(e_0)]^{-\frac{1}{2}}A(-1,e_0)[A(-1,e_0)+P_1(e_0)]^{-\frac{1}{2}}x_{\bb_0,e_0},x_{\bb_0,e_0}\>=0.
\ee
Let $y_0=[A(-1,e_0)+P_1(e_0)]^{-\frac{1}{2}}x_{\bb_0,e_0}$, and plugging it into (\ref{Bx.x=0}), we obtain
\be
\<A(-1,e_0)y_0,y_0\>=0,
\ee
and hence by Lemma 4.1, we must have
\be
y_0=(c_1(1+e_0\cos t),c_2(1+e_0\cos t))^T
\ee
for some constants $c_1,c_2\in\C$. Moreover, we have
\be
x_{\bb_0,e_0}=[A(-1,e_0)+P_1(e_0)]^{\frac{1}{2}}y_0=(c_1(1+e_0\cos t),c_2(1+e_0\cos t))^T=y_0.
\ee
Then $B(\bb_0,e_0)x_{\bb_0,e_0}=f(\bb_0,e_0)x_{\bb_0,e_0}$ reads
\bea
f(\bb_0,e_0)y_0&=&[A(-1,e_0)+P_1(e_0)]^{\frac{1}{2}}(f(\bb_0,e_0)x_{\bb_0,e_0}) \nonumber
\\
&=&[A(-1,e_0)+P_1(e_0)]^{\frac{1}{2}}B(\bb_0,e_0)x_{\bb_0,e_0}  \nonumber
\\
&=&\left(\frac{I_2+3S(t)}{2(1+e_0\cos t)}-\frac{P_1(e_0)}{\bb_0+1}\right)y_0  \nonumber
\\
&=&\frac{I_2+3S(t)}{2}\left(\matrix{c_1\cr c_2}\right)-\frac{1}{\bb_0+1}y_0,
\eea
this is impossible unless $c_1=c_2=0$.
Therefore $\frac{\partial f(\bb,e)}{\partial\bb}-\frac{1}{(\bb+1)^2}\ne0$,
and then apply the implicit function theorem to (\ref{f.be}),
$\bb_i(1,e)$ is real analytical functions of $e$.
\hb

Moreover, we have

\begin{theorem}\label{Th:orth}
Every $1$-degenerate curve must start from point $(\hat\bb_n,0),n\ge1$ and is orthogonal to the $\bb$-axis.
\end{theorem}

{\bf Proof.} Let $(\bb(e),e)$ be one of such curves
(i.e., one of $(\bb_i(1,e),e),i\in\N$. later, we will show that the two curves coincide)
which starts from $\bb(0)=\hat\bb_n$ with $e\in(-\epsilon,\epsilon)$ for some small $\epsilon>0$ and $x_e\in\bar{D}(1,2\pi)$
be the corresponding eigenvector, that is
\be
A(\bb(e),e)x_e=0.
\ee
Without loose of generality, by Remark 4.8, we suppose
\be
x_0=R(t)(a_n\sin nt,\cos nt)^T\lb{7.20a}
\ee
and
$$z=(a_n\sin nt,\cos nt)^T.$$
There holds
\be
\<A(\bb(e),e)x_e,x_e\>=0.\label{Axx1}
\ee

Differentiating both side of (\ref{Axx1}) with respect to $e$ yields
$$ \bb'(e)\<\frac{\pt}{\pt \bb}A(\bb(e),e)x_e,x_e\> + (\<\frac{\pt}{\pt e}A(\bb(e),e)x_e,x_e\>
       + 2\<A(\bb(e),e)x_e,x'_e\> = 0,  $$
where $\bb'(e)$ and $x'_e$ denote the derivatives with respect to $e$. Then evaluating both
sides at $e=0$ yields
\be  \bb'(0)\<\frac{\pt}{\pt \bb}A(\hat\bb_n,0)x_0,x_0\> + \<\frac{\pt}{\pt e}A(\hat\bb_n,0)x_0,x_0\> = 0. \lb{7.21a}\ee
Then by the definition (\ref{2.29}) of $A(\bb,e)$ we have
\bea
\left.\frac{\pt}{\pt\bb}A(\bb,e)\right|_{(\bb,e)=(\hat\bb_n,0)}
    &=& \left.R(t)\frac{\pt}{\pt\bb}K_{\bb,e}(t)\right|_{(\bb,e)=(\hat\bb_n,0)}R(t)^T,  \lb{7.22a}\\
\left.\frac{\pt}{\pt e}A(\bb,e)\right|_{(\bb,e)=(\hat\bb_n,0)}
    &=& \left.R(t)\frac{\pt}{\pt e}K_{\bb,e}(t)\right|_{(\bb,e)=(\hat\bb_n,0)}R(t)^T,  \lb{7.23a}\eea
where $R(t)$ is given in \S 2.1. By direct computations from the definition of $K_{\bb,e}(t)$ in
(\ref{2.20}), we obtain
\bea
&& \frac{\pt}{\pt\bb}K_{\bb,e}(t)\left|_{(\bb,e)=(\hat\bb_n,0)}
       = \left(\matrix{2 & 0\cr
                                          0 &  -1\cr}\right),\right.   \lb{7.24a}\\
&& \frac{\pt}{\pt e}K_{\bb,e}(t)\left|_{(\bb,e)=(\hat\bb_n,0)}
       = {-\cos t}\left(\matrix{2\hat\bb_n+3 & 0 \cr
                                        0 & -\hat\bb_n\cr}\right).\right.   \lb{7.25a}\eea
Therefore from (\ref{7.20a}) and (\ref{7.22a})-(\ref{7.25a}) we have
\bea  \<\frac{\pt}{\pt\bb}A(\hat\bb_n,0)x_0,x_0\>
&=& \<\frac{\pt}{\pt\bb}K_{\hat\bb_n,0}z,z\>    \nn\\
&=& \int_0^{2\pi}[2a_n^2\sin^2nt
                -\cos^2nt]dt  \nn\\
&=& \pi(2a_n^2-1),  \lb{7.26a}\eea
and
\bea  \<\frac{\pt}{\pt e}A(\hat\bb_n,0)x_0,x_0\>
&=& \<\frac{\pt}{\pt e}K_{\hat\bb_n,0}z,z\>   \nn\\
&=& \int_0^{2\pi}[-(2\bb_n+3)a_n^2\cos{t}\sin^2nt
                +\hat\bb_n\cos{t}\cos^2nt]dt  \nn\\
&=& 0.  \lb{7.27a}\eea
Therefore by (\ref{7.21a}) and (\ref{7.26a})-(\ref{7.27a}), together with $a_n^2\ne1/2$ which from Remark 4.8,
we obtain
\be  \beta'(0) = 0.  \lb{7.28a}\ee
Thus the theorem is proved.\hb

\subsection{The $\omega=-1$ degenerate curves}

Recall $A(-1,e)$ is positive definite on $\overline{D}(\om,2\pi)$ for $\om\ne1$.
Now we set
\be\lb{tilde.B}
\tilde{B}(e,\om)=A(-1,e)^{-\frac{1}{2}}
         \frac{I_2+3S(t)}{2(1+e\cos t)}
         A(-1,e)^{-\frac{1}{2}}.
\ee
Then we have
\begin{lemma}\label{L4.13}
For $-1< e<1$, $A(\bb,e)$ is $\om$-degenerate if and only if $-\frac{1}{\bb+1}$ is an eigenvalue of $\tilde{B}(e,\om)$.
\end{lemma}
{\bf Proof.} Suppose $A(\bb,e)x=0$ holds for some $x\in\overline{D}(\om,2\pi)$.
Let $y=A(-1,e)^{\frac{1}{2}}x$. Then by (\ref{tilde.B}) we obtain
\bea  A(-1,e)^{\frac{1}{2}}\left(\frac{1}{\bb+1}+\tilde{B}(e,\om)\right)y(t)
&=& \left(\frac{A(-1,e)}{\bb+1}+\frac{I_2+3S(t)}{2(1+e\cos t)}\right)x(t) \nonumber\\
&=& \frac{1}{\bb+1}A(\bb,e)x  \nonumber\\
&=& 0.    \label{A.tilde.B}\eea

Conversely, if $(\frac{1}{\bb+1}+\tilde{B}(e,\om))y=0$,
then $x=A(-1,e)^{-\frac{1}{2}}y$ is an eigenfunction of $A(\bb,e)$ belonging to the eigenvalue $0$
by our computations (\ref{A.tilde.B}).
\hb

For convenience, we define
\be
\bb_0(1,e)\equiv0 \quad\forall e\in[0,1).
\ee
We first have
\begin{theorem}\label{Th:analytic.om.orign}
For $\om\ne1$, there exists two analytic $\om$-degenerate curves $(h_i(e),e)$ in $e\in(-1,1)$ with $i=1,2$
such that $\bb_{2n}(1,e)<h_i(e)<\bb_{2n+1}(1,e),n\ge0$.
Specially, each $h_i(e)$ is a real analytic function in $e\in(-1,1)$ and $\bb_{2n}(1,e)<h_i(e)<\bb_{2n+1}(1,e)$.
Moreover, $\ga_{h_i(e),e}(2\pi)$ is $\om$-degenerate for $\om\in\U\backslash\{1\}$ and $i=1,2$.
\end{theorem}

{\bf Proof.}For $\bb\in(\bb_{2n}(1,e),\bb_{2n+1}(1,e))$, from Theorem 1.5 (ix)-(xiv),
we have
\be
i_1(\ga_{\bb,e})=2n+3,\quad \nu_1(\ga_{\bb,e})=0.
\ee
Moreover, from Theorem 1.5 (viii), we have
\be\label{boundary}
\ga_{\bb,e}\approx I_2\diamond D(2),\quad \bb=\bb_{2n}(1,e)\;\;{\rm or}\;\;\bb_{2n+1}(1,e).
\ee
Then for $\om\in\U\backslash\{1\}$, we have
\bea
i_\om(\ga_{\bb_{2n}(1,e),e})&=&i_1(\ga_{\bb_{2n}(1,e),e})+S_{\ga_{\bb_{2n}(1,e),e}(2\pi)}^+(1)
\nn\\
&=&2n+1+S_{I_2}^+(1)
\nn\\
&=&2n+2. \label{left.om.index}
\eea
Similarly, we have
\be
i_\om(\ga_{\bb_{2n+1}(1,e),e})=2n+4.\label{right.om.index}
\ee
Therefore, by Lemma \ref{Lemma:increasing.of.index}, it shows that, for fixed $e\in(-1,1)$,
there are exactly two values $\bb=h_1(e)$ and $h_2(e)$ in the interval $[\bb_{2n}(1,e),\bb_{2n+1}(1,e)]$
at which (\ref{A.tilde.B}) is satisfied, and then $\bar{A}(\bb,e)$ at these two values is $\om$-degenerate.
Note that these two $\bb$ values are possibly equal to each other at some $e$.
Moreover, (\ref{boundary}) implies that $h_i(e)\ne \bb_{2n}(1,e)$ and $\bb_{2n+1}(1,e)$ for $i=1, 2$.

By Lemma \ref{L4.13}, $-\frac{1}{\bb_i(\om,e)+1}$ is an eigenvalue of $\tilde{B}(e,\om)$.
Note that $\tilde{B}(e,\om)$ is a compact operator and self adjoint when $e$ are real.
Moreover, it depends analytically on $e$. By \cite{Ka}(Theorem 3.9 in p.392), we know that
$-\frac{1}{\bb_i(\om,e)+1}$ is analytic in $e$ for each $i\in\N$. This in turn implies that
both $h_1(e)$ and $h_2(e)$ are real analytic functions in $e$. \hb

By the definition of $\bb_n(\om,e)$ in (\ref{bb_s}), together with (\ref{index.of.0e}),
(\ref{null.index.of.0e}), (\ref{left.om.index}) and (\ref{right.om.index}),
we have
\bea
\bb_{2n+1}(\om,e)=\min\{h_1(e),h_2(e)\},
\\
\bb_{2n+2}(\om,e)=\max\{h_1(e),h_2(e)\}.
\eea
Thus we have the following theorem:

\begin{theorem}\label{Th:analytic.om}
For $\om\ne 1$, every $\om$-degenerate curve $(\bb_n(\om,e),e)$ in $e\in (-1,1)$ is a piecewise
analytic function.
\end{theorem}

For $\hat{\bb}_{n+\frac{1}{2}}$ defined by (\ref{om_n}),
$A(\hat\bb_{n+\frac{1}{2}},0)$ is degenerate and by (\ref{null.-1.index.of.b0}),
$\dim\ker A(\hat\bb_{n+\frac{1}{2}},0)=v_{-1}(\ga_{\hat\bb_{n+\frac{1}{2}},0})=2$.
By the definition of (\ref{2.11}), we have
$R(t)\left(\matrix{\tilde{a}_n\sin (n+\frac{1}{2})t\cr \cos (n+\frac{1}{2})t}\right)\in\ol{D}(-1,2\pi)$
for any constant $\tilde{a}_n$.

Moreover,
$A(\bb,0)R(t)\left(\matrix{\tilde{a}_n\sin (n+\frac{1}{2})t\cr \cos (n+\frac{1}{2})t}\right)=0$ reads
\be
\left\{
\begin{array}{cr}
(n+\frac{1}{2})^2\tilde{a}_n-2(n+\frac{1}{2})+(2\bb+3)\tilde{a}_n&=0,
\\
(n+\frac{1}{2})^2-2(n+\frac{1}{2})\tilde{a}_n-\bb&=0.
\end{array}
\right.
\ee
Then $2\bb^2-((n+\frac{1}{2})^2-3)\bb-(n+\frac{1}{2})^2((n+\frac{1}{2})^2-1)=0$ which yields $\bb=\hat\bb_{n+\frac{1}{2}}$ again and
\be\lb{tilde.a}
\tilde{a}_n=\frac{(n+\frac{1}{2})^2-\hat\bb_{n+\frac{1}{2}}}{2n+1}.
\ee
Then we have
$R(t)\left(\matrix{\tilde{a}_n\sin (n+\frac{1}{2})t\cr \cos (n+\frac{1}{2})t}\right)\in\ker A(\hat\bb_{n+\frac{1}{2}},0)$.
Similarly
$R(t)\left(\matrix{\tilde{a}_n\cos (n+\frac{1}{2})t\cr -\sin (n+\frac{1}{2})t}\right)\in\ker A(\hat\bb_{n+\frac{1}{2}},0)$,
therefore we have
\begin{equation}
\ker A(\hat\bb_{n+\frac{1}{2}},0)=\span\left\{
R(t)\left(\matrix{\tilde{a}_n\sin (n+\frac{1}{2})t\cr \cos (n+\frac{1}{2})t}\right),\quad
R(t)\left(\matrix{\tilde{a}_n\cos (n+\frac{1}{2})t\cr -\sin (n+\frac{1}{2})t}\right)
\right\}.
\lb{ker.A.of.-1}
\end{equation}
Indeed, we have the following theorem:

\begin{theorem}\label{Th:orth.om}
Every $-1$-degenerate curve must start from the point $(\hat\bb_{n+\frac{1}{2}},0),n\ge1$ and is
orthogonal to the $\bb$-axis.
\end{theorem}

{\bf Proof.} Similarly to Lemma \ref{L4.9}, we have
\begin{equation}
\bb_n(-1,0)=
\left\{
\begin{array}{l}
\hat\bb_{m+\frac{1}{2}}, \quad {\rm if}\;\;n = 2m-1,
\\
\hat\bb_{m+\frac{1}{2}}, \quad {\rm if}\;\;n = 2m.
\end{array}
\right.
\lb{bb_n0.of.-1}
\end{equation}
Thus every $-1$-degenerate curve $(\bb(-1,e),e)$ must start from point $(\hat\bb_{n+\frac{1}{2}},0)$.

Now let $(\bb(e),e)$ be one of such curves
(i.e., one of $(\bb_i(-1,e),e),i\in\N$.)
which starts from $\bb(0)=\hat\bb_{n+\frac{1}{2}}$ with $e\in(-\epsilon,\epsilon)$ for some small $\epsilon>0$ and $x_e\in\bar{D}(1,2\pi)$
be the corresponding eigenvector, that is
\be
A(\bb(e),e)x_e=0.
\ee
Without loose of generality, by (\ref{ker.A.of.-1}), we suppose
$$
z = (\tilde{a}_n\sin (n+\frac{1}{2})t,\cos (n+\frac{1}{2})t)^T
$$
and
\be
x_0=R(t)z=R(t)(\tilde{a}_n\sin (n+\frac{1}{2})t,\cos (n+\frac{1}{2})t)^T.\lb{4.71}
\ee
There holds
\be
\<A(\bb(e),e)x_e,x_e\>=0.\label{Axx-1}
\ee

Differentiating both side of (\ref{Axx-1}) with respect to $e$ yields
$$ \bb'(e)\<\frac{\pt}{\pt \bb}A(\bb(e),e)x_e,x_e\> + (\<\frac{\pt}{\pt e}A(\bb(e),e)x_e,x_e\>
       + 2\<A(\bb(e),e)x_e,x'_e\> = 0,  $$
where $\bb'(e)$ and $x'_e$ denote the derivatives with respect to $e$. Then evaluating both
sides at $e=0$ yields
\be  \bb'(0)\<\frac{\pt}{\pt \bb}A(\hat\bb_{n+\frac{1}{2}},0)x_0,x_0\>
   + \<\frac{\pt}{\pt e}A(\hat\bb_{n+\frac{1}{2}},0)x_0,x_0\> = 0.
\lb{4.73}
\ee
Then by the definition (\ref{2.29}) of $A(\bb,e)$ we have
\bea
\left.\frac{\pt}{\pt\bb}A(\bb,e)\right|_{(\bb,e)=(\hat\bb_{n+\frac{1}{2}},0)}
    &=& \left.R(t)\frac{\pt}{\pt\bb}K_{\bb,e}(t)\right|_{(\bb,e)=(\hat\bb_{n+\frac{1}{2}},0)}R(t)^T,  \lb{4.74}\\
\left.\frac{\pt}{\pt e}A(\bb,e)\right|_{(\bb,e)=(\hat\bb_{n+\frac{1}{2}},0)}
    &=& \left.R(t)\frac{\pt}{\pt e}K_{\bb,e}(t)\right|_{(\bb,e)=(\hat\bb_{n+\frac{1}{2}},0)}R(t)^T,  \lb{4.75}\eea
where $R(t)$ is given in \S 2.1. By direct computations from the definition of $K_{\bb,e}(t)$ in
(\ref{2.20}), we obtain
\bea
&& \frac{\pt}{\pt\bb}K_{\bb,e}(t)\left|_{(\bb,e)=(\hat\bb_{n+\frac{1}{2}},0)}
       = \left(\matrix{2 & 0\cr
                                          0 &  -1\cr}\right),\right.   \lb{4.76}\\
&& \frac{\pt}{\pt e}K_{\bb,e}(t)\left|_{(\bb,e)=(\hat\bb_{n+\frac{1}{2}},0)}
       = {-\cos t}\left(\matrix{2\hat\bb_{n+\frac{1}{2}}+3 & 0 \cr
                                        0 & -\hat\bb_{n+\frac{1}{2}}\cr}\right).\right.   \lb{4.77}\eea
Therefore from (\ref{4.71}) and (\ref{4.74})-(\ref{4.77}) we have
\bea  \<\frac{\pt}{\pt\bb}A(\hat\bb_{n+\frac{1}{2}},0)x_0,x_0\>
&=& \<\frac{\pt}{\pt\bb}K_{\hat\bb_{n+\frac{1}{2}},0}z,z\>    \nn\\
&=& \int_0^{2\pi}[2\tilde{a}_n^2\sin^2(n+\frac{1}{2})t
                -\cos^2(n+\frac{1}{2})t]dt  \nn\\
&=& \pi(2\tilde{a}_n^2-1),  \lb{4.78}\eea
and for $n\ge1$,
\bea  \<\frac{\pt}{\pt e}A(\hat\bb_{n+\frac{1}{2}},0)x_0,x_0\>
&=& \<\frac{\pt}{\pt e}K_{\hat\bb_{n+\frac{1}{2}},0}z,z\>   \nn\\
&=& \int_0^{2\pi}[-(2\bb_{n+\frac{1}{2}}+3)\tilde{a}_n^2\cos{t}\sin^2({n+\frac{1}{2}})t
                +\hat\bb_{n+\frac{1}{2}}\cos{t}\cos^2({n+\frac{1}{2}})t]dt  \nn\\
&=& 0.  \lb{4.79}\eea
Therefore by (\ref{4.73}) and (\ref{4.78})-(\ref{4.79}),
together with $\tilde{a}_n^2\ne1/2$ which from (\ref{om_n.5}) and (\ref{tilde.a}),
we obtain
\be  \beta'(0) = 0.  \lb{4.80}\ee
Thus the theorem is proved.\hb

\setcounter{equation}{0}
\section{Appendix}

\subsection{On $\dl$ and $\bb$.}

\begin{lemma}
Let $(m_1,m_2,m_3)\in\R^3$ satisfying (\ref{nomorlize.the.masses}),
and $x$ be any solution of the quintic polynomial (\ref{quintic.polynomial}),
then there holds
\begin{eqnarray}
&&\frac{m_3(1+x)^3(m_2+m_1+m_1x)^2+m_1x^3(1+x)^3(m_3+m_3x+m_2x)^2+m_2x^3(m_1x-m_3)^2}
        {x^2(1+x)^2[m_2m_3x^2+(m_1m_2+m_2m_3+m_3m_1)x+m_1m_2]}
\nonumber\\
&&=1+\frac{m_1(3x^2+3x+1)+m_3x^2(x^2+3x+3)}
          {x^2+m_2[(x+1)^2(x^2+1)-x^2]}.\label{beta.delta}
\end{eqnarray}
\end{lemma}

{\bf Proof. }Firstly let's define
\bea
q_0&=&(x+1)^2(x^2+1)-x^2=x^4+2x^3+x^2+2x+1,\label{q_0}
\\
r&=&\frac{x^3(x^2+3x+3)}{(x+1)q_0}=\frac{x^3(x^2+3x+3)}{(x+1)(x^4+2x^3+x^2+2x+1)},\label{r}
\\
p_1&=&m_3(1+x)^3(m_2+m_1+m_1x)^2+m_1x^3(1+x)^3(m_3+m_3x+m_2x)^2+m_2x^3(m_1x-m_3)^2,
\\
q_1&=&x^2(1+x)^2[m_2m_3x^2+(m_1m_2+m_2m_3+m_3m_1)x+m_1m_2],
\\
p_2&=&m_1(3x^2+3x+1)+m_3x^2(x^2+3x+3)+x^2+m_2[(x+1)^2(x^2+1)-x^2],
\\
q_2&=&x^2+m_2[(x+1)^2(x^2+1)-x^2].\label{q_2}
\eea
By (\ref{nomorlize.the.masses})and (\ref{quintic.polynomial}), we can represent $m_2$ and $m_3$ by $m_1$ and $x$:
\bea
m_1&=&=-\frac{1}{x+1}m_2+r,
\\
m_3&=&-\frac{x}{x+1}m_2+1-r.
\eea
Therefore, we use $m_2$ and $x$ as our parameters. Moreover, by (\ref{r}), we have
\bea
1-r&=&1-\frac{x^3(x^2+3x+3)}{(x+1)(x^4+2x^3+x^2+2x+1)}=\frac{3x^2+3x+1}{(x+1)(x^4+2x^3+x^2+2x+1)},\label{one.minus.r}
\\
m_1x-m_3&=&x(-\frac{1}{x+1}m_2+r)-(-\frac{x}{x+1}m_2+1-r)=(1+x)r-1.
\eea
Using (\ref{q_0})-(\ref{one.minus.r}), by directly computation, we have
\bea
q_2&=&x^2+m_2q_0=q_0(m_2+\frac{x^2}{q_0}),\label{q_2'}
\\
p_2&=&(-\frac{1}{x+1}m_2+r)(3x^2+3x+1)+(-\frac{x}{x+1}m_2+1-r)x^2(x^2+3x+3)+q_0(m_2+\frac{x^2}{q_0})
\nonumber\\
&=&-\frac{x^5+3x^4+3x^3+3x^2+3x+1}{x+1}m_2+\frac{(x^2+3x+3)(3x^2+3x+1)(x^3+x^2)}{(x+1)q_0}+q_0(m_2+\frac{x^2}{q_0})
\nonumber\\
&=&-q_0m_2+\frac{x^2(x^2+3x+3)(3x^2+3x+1)}{q_0}+q_0(m_2+\frac{x^2}{q_0})
\nonumber\\
&=&\frac{x^2(x^2+3x+3)(3x^2+3x+1)}{q_0}+x^2
\nonumber\\
&=&\frac{2x^2}{q_0}(2x^4+7x^3+10x^2+7x+2)
\nonumber\\
&=&\frac{2x^2(x+1)^2(2x^2+3x+2)}{q_0},
\\
q_1&=&x^2(x+1)^2[m_2m_3x(x+1)+m_1m_2(x+1)+m_1m_3x]
\nonumber\\
&=&x^2(x+1)^2[m_2x(-m_2x+\frac{3x^2+3x+1}{q_0})+m_2(-m_2+\frac{x^3(x^2+3x+3)}{q_0})
\nonumber\\
&&+(-\frac{1}{x+1}m_2+r)(-\frac{x}{x+1}m_2+1-r)x]
\nonumber\\
&=&x^2(x+1)^2[-\frac{q_0}{(x+1)^2}m_2^2+\frac{x^2(2x^4+10x^3+18x^2+10x+2)}{(x+1)^2q_0}m_2+\frac{x^4(x^2+3x+3)(3x^2+3x+1)}{(x+1)^2q_0^2}]
\nonumber\\
&=&-x^2q_0[m_2^2-\frac{x^2(2x^4+10x^3+18x^2+10x+2)}{q_0^2}m_2-\frac{x^4(x^2+3x+3)(3x^2+3x+1)}{q_0^3}]
\nonumber\\
&=&-x^2q_0[m_2^2+(\frac{x^2}{q_0}-\frac{x^2(2x^4+10x^3+18x^2+10x+2+q_0)}{q_0^2})m_2-\frac{x^4(x^2+3x+3)(3x^2+3x+1)}{q_0^3}]
\nonumber\\
&=&-x^2q_0[m_2^2+(\frac{x^2}{q_0}-\frac{x^2(x^2+3x+3)(3x^2+3x+1)}{q_0^2})m_2-\frac{x^4(x^2+3x+3)(3x^2+3x+1)}{q_0^3}]
\nonumber\\
&=&-x^2q_0(m_2+\frac{x^2}{q_0})[m_2-\frac{x^2(x^2+3x+3)(3x^2+3x+1)}{q_0^2}],
\\
p_1&=&m_3(x+1)^3(1-m_3+m_1x)+m_1x^3(x+1)^3(m_3+x-m_1x)^2+m_2x^3(m_1x-m_3)^2
\nonumber\\
&=&m_3(x+1)^3[(x+1)r]^2+m_1x^3(x+1)^3[(x+1)(1-r)]^2+m_2x^3[(x+1)r-1]^2
\nonumber\\
&=&(-\frac{x}{x+1}m_2+1-r)(x+1)^5r^2+(-\frac{1}{x+1}m_2+r)x^3(x+1)^5(1-r)^2+m_2x^3[(x+1)r-1]^2
\nonumber\\
&=&-[x(x+1)^4r^2+x^3(x+1)^4(1-r)^2-x^3(xr+r-1)^2]m_2+[(x+1)^5r^2(1-r)+x^3(x+1)^5r(1-r)^2]
\nonumber\\
&=&-\frac{x^3(x+1)^2}{q_0^2}[x^4(x^2+3x+3)^2+(3x^2+3x+1)^2-(x+1)^2(x^3-1)^2]m_2
\nonumber\\
&&+(x+1)^5r(1-r)[r+x^3(1-r)]
\nonumber\\
&=&-\frac{x^3(x+1)^2}{q_0^2}[2x(2^6+7x^5+10x^4+11x^3+10x^2+7x+2)]m_2
\nonumber\\
&&+\frac{x^3(x+1)^3(x^2+3x+3)(3x^2+3x+1)}{q_0^2}\frac{2x^3(2x^2+3x+2)}{(x+1)q_0}
\nonumber\\
&=&-\frac{x^3(x+1)^2}{q_0^2}[2xq_0(2x^2+3x+2)]m_2+\frac{x^6(x+1)^2(x^2+3x+3)(3x^2+3x+1)(2x^2+3x+2)}{q_0^3}
\nonumber\\
&=&-\frac{2x^4(x+1)^2(2x^2+3x+2)}{q_0}[m_2-\frac{x^2(x^2+3x+3)(3x^2+3x+1)}{q_0^2}].\label{p_1'}
\eea
Thus by (\ref{q_2'})-(\ref{p_1'}), (\ref{beta.delta}) holds.\hb

\medskip

\subsection{$\omega$-Maslov-type indices and $\omega$-Morse indices}

Let $(\R^{2n},\Omega)$ be the standard symplectic vector space with coordinates
$(x_1,...,x_n,y_1,...,y_n)$ and the symplectic form $\Omega=\sum_{i=1}^{n}dx_i \wedge dy_i$.
Let $J=\left(\matrix{0&-I_n\cr
                 I_n&0\cr}\right)$ be the standard symplectic matrix, where $I_n$
is the identity matrix on $\R^n$.

As usual, the symplectic group $\Sp(2n)$ is defined by
$$ \Sp(2n) = \{M\in {\rm GL}(2n,\R)\,|\,M^TJM=J\}, $$
whose topology is induced from that of $\R^{4n^2}$. For $\tau>0$
we are interested in paths in $\Sp(2n)$:
$$ \P_{\tau}(2n) = \{\ga\in C([0,\tau],\Sp(2n))\,|\,\ga(0)=I_{2n}\}, $$
which is equipped with the topology induced from that of $\Sp(2n)$.
For any $\om\in\U$ and $M\in\Sp(2n)$, the following real function was
introduced in \cite{Lon2}:
$$ D_{\om}(M) = (-1)^{n-1}\ol{\om}^n\det(M-\om I_{2n}). $$
Thus for any $\om\in\U$ the following codimension $1$ hypersurface
in $\Sp(2n)$ is defined (\cite{Lon2}):
$$ \Sp(2n)_{\om}^0 = \{M\in\Sp(2n)\,|\, D_{\om}(M)=0\}.  $$
For any $M\in \Sp(2n)_{\om}^0$, we define a co-orientation of
$\Sp(2n)_{\om}^0$ at $M$ by the positive direction
$\frac{d}{dt}Me^{t J}|_{t=0}$ of the path $Me^{t J}$ with $0\le t\le
\varepsilon$ and $\varepsilon$ being a small enough positive number. Let
\bea
\Sp(2n)_{\om}^{\ast} &=& \Sp(2n)\bs \Sp(2n)_{\om}^0,   \nn\\
\P_{\tau,\om}^{\ast}(2n) &=&
      \{\ga\in\P_{\tau}(2n)\,|\,\ga(\tau)\in\Sp(2n)_{\om}^{\ast}\}, \nn\\
\P_{\tau,\om}^0(2n) &=& \P_{\tau}(2n)\bs \P_{\tau,\om}^{\ast}(2n). \nn\eea
For any two continuous paths $\xi$ and $\eta:[0,\tau]\to\Sp(2n)$ with
$\xi(\tau)=\eta(0)$, we define their concatenation by:
$$ \eta\ast\xi(t) = \left\{\matrix{
            \xi(2t), & \quad {\rm if}\;0\le t\le \tau/2, \cr
            \eta(2t-\tau), & \quad {\rm if}\; \tau/2\le t\le \tau. \cr}\right. $$

As in \cite{Lon4}, for $\lm\in\R\bs\{0\}$, $a\in\R$, $\th\in (0,\pi)\cup (\pi,2\pi)$,
$b=\left(\matrix{b_1 & b_2\cr
                 b_3 & b_4\cr}\right)$ with $b_i\in\R$ for $i=1, \ldots, 4$, and $c_j\in\R$
for $j=1, 2$, we denote respectively some normal forms by
\bea
&& D(\lm)=\left(\matrix{\lm & 0\cr
                         0  & \lm^{-1}\cr}\right), \qquad
   R(\th)=\left(\matrix{\cos\th & -\sin\th\cr
                        \sin\th  & \cos\th\cr}\right),  \nn\\
&& N_1(\lm, a)=\left(\matrix{\lm & a\cr
                             0   & \lm\cr}\right), \qquad
   N_2(e^{\sqrt{-1}\th},b) = \left(\matrix{R(\th) & b\cr
                                           0      & R(\th)\cr}\right),  \nn\\
&& M_2(\lm,c)=\left(\matrix{\lm &   1 &       c_1 &         0 \cr
                              0 & \lm &       c_2 & (-\lm)c_2 \cr
                              0 &   0 &  \lm^{-1} &         0 \cr
                              0 &   0 & -\lm^{-2} &  \lm^{-1} \cr}\right). \nn\eea
Here $N_2(e^{\sqrt{-1}\th},b)$ is {\bf trivial} if $(b_2-b_3)\sin\th>0$, or {\bf non-trivial}
if $(b_2-b_3)\sin\th<0$, in the sense of Definition 1.8.11 on p.41 of \cite{Lon4}. Note that
by Theorem 1.5.1 on pp.24-25 and (1.4.7)-(1.4.8) on p.18 of \cite{Lon4}, when $\lm=-1$ there hold
\bea
c_2 \not= 0 &{\rm if\;and\;only\;if}\;& \dim\ker(M_2(-1,c)+I)=1, \nn\\
c_2 = 0 &{\rm if\;and\;only\;if}\;& \dim\ker(M_2(-1,c)+I)=2. \nn\eea
Note that we have $N_1(\lm,a)\approx N_1(\lm, a/|a|)$ for $a\in\R\bs\{0\}$ by symplectic coordinate
change, because
$$ \left(\matrix{1/\sqrt{|a|} & 0\cr
                           0  & \sqrt{|a|}\cr}\right)
   \left(\matrix{\lm & a\cr
                  0  & \lm\cr}\right)
   \left(\matrix{\sqrt{|a|} & 0\cr
                           0  & 1/\sqrt{|a|}\cr}\right) = \left(\matrix{\lm & a/|a|\cr
                                                                         0  & \lm\cr}\right). $$

\begin{definition}\lb{D2.1} (\cite{Lon2}, \cite{Lon4})
For any $\om\in\U$ and $M\in \Sp(2n)$, define
\be \nu_{\om}(M)=\dim_{\C}\ker_{\C}(M - \om I_{2n}).  \lb{2.2}\ee

For every $M\in \Sp(2n)$ and $\om\in\U$, as in Definition 1.8.5 on p.38 of \cite{Lon4}, we define the
{\bf $\om$-homotopy set} $\Om_{\om}(M)$ of $M$ in $\Sp(2n)$ by
$$  \Om_{\om}(M)=\{N\in\Sp(2n)\,|\, \nu_{\om}(N)=\nu_{\om}(M)\},  $$
and the {\bf homotopy set} $\Om(M)$ of $M$ in $\Sp(2n)$ by
\bea  \Om(M)=\{N\in\Sp(2n)\,&|&\,\sg(N)\cap\U=\sg(M)\cap\U,\,{\it and}\; \nn\\
         &&\qquad \nu_{\lm}(N)=\nu_{\lm}(M)\qquad\forall\,\lm\in\sg(M)\cap\U\}.  \nn\eea
We denote by $\Om^0(M)$ (or $\Om^0_{\om}(M)$) the path connected component of $\Om(M)$ ($\Om_{\om}(M)$)
which contains $M$, and call it the {\bf homotopy component} (or $\om$-{\bf homtopy component}) of $M$ in
$\Sp(2n)$. Following Definition 5.0.1 on p.111 of \cite{Lon4}, for $\om\in \U$ and $\ga_i\in \P_{\tau}(2n)$
with $i=0, 1$, we write $\ga_0\sim_{\om}\ga_1$ if $\ga_0$ is homotopic to $\ga_1$ via
a homotopy map $h\in C([0,1]\times [0,\tau], \Sp(2n))$ such that $h(0)=\ga_0$, $h(1)=\ga_1$, $h(s)(0)=I$,
and $h(s)(\tau)\in \Om_{\om}^0(\ga_0(\tau))$ for all $s\in [0,1]$. We write also $\ga_0\sim \ga_1$, if
$h(s)(\tau)\in \Om^0(\ga_0(\tau))$ for all $s\in [0,1]$ is further satisfied. We write $M\approx N$, if
$N\in \Om^0(M)$.
\end{definition}

Following Definition 1.8.9 on p.41 of \cite{Lon4}, we call the above matrices $D(\lm)$, $R(\th)$, $N_1(\lm,a)$
and $N_2(\om,b)$ basic normal forms of symplectic matrices. As proved in \cite{Lon2} and \cite{Lon3} (cf.
Theorem 1.9.3 on p.46 of \cite{Lon4}), every $M\in\Sp(2n)$ has its basic normal form decomposition in $\Om^0(M)$
as a $\dm$-sum of these basic normal forms. Here the $\dm$-sum is introduced in the above Section 1. This is
very important when we derive basic normal forms for $\ga_{\bb,e}(2\pi)$ to compute the $\om$-index
$i_{\om}(\ga_{\bb,e})$ of the path $\ga_{\bb,e}$ later in this paper.

We define a special continuous symplectic path $\xi_n\subset {\Sp}(2n)$ by
\be \xi_n(t) = \left(\matrix{2-\frac{t}{\tau} & 0 \cr
                             0 &  (2-\frac{t}{\tau})^{-1}\cr}
               \right)^{\dm n} \qquad {\rm for}\;0\le t\le \tau.  \lb{2.3}\ee

\begin{definition} (\cite{Lon2}, \cite{Lon4})\lb{D2.2}
{For any $\tau>0$ and $\ga\in \P_{\tau}(2n)$, define
\be \nu_{\om}(\ga)= \nu_{\om}(\ga(\tau)).  \lb{2.4}\ee

If $\ga\in\P_{\tau,\om}^{\ast}(2n)$, define
\be i_{\om}(\ga) = [\Sp(2n)_{\om}^0: \ga\ast\xi_n],  \lb{2.5}\ee
where the right hand side of (\ref{2.5}) is the usual homotopy intersection number, and
the orientation of $\ga\ast\xi_n$ is its positive time direction under homotopy with
fixed end points.

If $\ga\in\P_{\tau,\om}^0(2n)$, we let $\mathcal{F}(\ga)$ be the set of all open
neighborhoods of $\ga$ in $\P_{\tau}(2n)$, and define
\be i_{\om}(\ga) = \sup_{U\in\mathcal{F}(\ga)}\inf\{i_{\om}(\beta)\,|\,
                       \beta\in U\cap\P_{\tau,\om}^{\ast}(2n)\}.      \lb{2.6}\ee
Then
$$ (i_{\om}(\ga), \nu_{\om}(\ga)) \in \Z\times \{0,1,\ldots,2n\}, $$
is called the index function of $\ga$ at $\om$. }
\end{definition}

\begin{definition} (\cite{Lon2}, \cite{Lon4})\lb{D2.3}
For any $M\in\Sp(2n)$ and $\om\in\U$, choosing $\tau>0$ and $\ga\in\P_\tau(2n)$ with $\ga(\tau)=M$,
we define
\be
S_M^{\pm}(\om)=\lim_{\epsilon\rightarrow0^+}\;i_{\exp(\pm\epsilon\sqrt{-1}\om)}(\ga)-i_\om(\ga).
\ee
They are called the splitting numbers of $M$ at $\om$.
\end{definition}
The splitting numbers $S_M^{\pm}(\om)$ measures the jumps between $i_\om(\ga)$
and $i_\lambda(\ga)$ with $\lambda\in\U$ near $\om$ from two sides of $\om$ in $\U$.
Therefore for any $\om_0=e^{\sqrt{-1}\th_0}\in\U$ with $0\le\th_0<2\pi$,
we denote by $\om_j$ with $1\le j\le p_0$ the eigenvalues of $M$ on $\U$ which are
distributed counterclockwise from $1$ to $\om_0$ and located strictly between $1$ and $\om_0$.
Then we have
\be i_{\om_0}(\ga)=i_1(\ga)+S_M^+(1)+\sum_{j=1}^{p_0}(-S_M^-(\om_j)+S_M^+(\om_j))-S_M^-(\om_0).  \lb{5.A1}\ee

\begin{lemma}
(Long, \cite{Lon4},p.198)
The integer valued splitting number pair $(S_M^+(\om),S_M^-(\om))$ defined for all
$(\om,M)\in\U\times\cup_{n\ge1}\Sp(2n)$ are uniquely determined by the following axioms:

$1^{\circ}$ (Homotopy invariant) $S_M^{\pm}(\om)=S_N^{\pm}(\om)$ for all $N\in\Omega^0(M)$.

$2^{\circ}$ (Symplectic additivity) $S_{M_1\diamond M_2}^{\pm}(\om)=S_{M_1}^{\pm}(\om)+S_{M_2}^{\pm}(\om)$
for all $M_i\in\Sp(2n_i)$ with $i=1$ and $2$.

$3^{\circ}$ (Vanishing) $S_M^{\pm}(\om)=0$ if $\om\not\in\sigma(M)$.

$4^{\circ}$ (Normality) $(S_M^+(\om),S_M^-(\om))$ coincides with the ultimate type of
$\om$ for $M$ when $M$ is any basic normal form.
\end{lemma}
Moreover, for $\om\in\C$ and $M\in\Sp(2n)$, we have
\be
S_M^+(\om)=S_M^-(\overline\om).
\ee

For the reader's convenience, we list the splitting numbers blow for all basic normal forms:

$\langle$1$\rangle$ $(S_M^+(1),S_M^-(1))=(1,1)$ for $M=N_1(1,b)$ with $b=1$ or $0$.

$\langle$2$\rangle$ $(S_M^+(1),S_M^-(1))=(0,0)$ for $M=N_1(1,-1)$.

$\langle$3$\rangle$ $(S_M^+(-1),S_M^-(-1))=(1,1)$ for $M=N_1(-1,b)$ with $b=-1$ or $0$.

$\langle$4$\rangle$ $(S_M^+(-1),S_M^-(-1))=(0,0)$ for $M=N_1(-1,1)$.

$\langle$5$\rangle$ $(S_M^+(e^{\sqrt{-1}\th}),S_M^-(e^{\sqrt{-1}\th}))=(0,1)$
for $M=R(\theta)$ with $\th\in(0,\pi)\cup(\pi2\pi)$.

$\langle$6$\rangle$ $(S_M^+(\om),S_M^-(\om))=(1,1)$ for $M=N_2(\om,b)$
being non-trivial with $\om=e^{\sqrt{-1}\th}\in\U\backslash\R$.

$\langle$7$\rangle$ $(S_M^+(\om),S_M^-(\om))=(0,0)$ for $M=N_2(\om,b)$
being trivial with $\om=e^{\sqrt{-1}\th}\in\U\backslash\R$.

$\langle$8$\rangle$ $(S_M^+(\om),S_M^-(\om))=(0,0)$ for $\om\in\U$ and $M=\in\Sp(2n)$
satisfying $\sigma(M)\cap\U=\emptyset$.

We refer to \cite{Lon4} for more details on this index theory of symplectic matrix paths
and periodic solutions of Hamiltonian system.

For $T>0$, suppose $x$ is a critical point of the functional
$$ F(x)=\int_0^TL(t,x,\dot{x})dt,  \qquad \forall\,\, x\in W^{1,2}(\R/T\Z,\R^n), $$
where $L\in C^2((\R/T\Z)\times \R^{2n},\R)$ and satisfies the
Legendrian convexity condition $L_{p,p}(t,x,p)>0$. It is well known
that $x$ satisfies the corresponding Euler-Lagrangian
equation:
\bea
&& \frac{d}{dt}L_p(t,x,\dot{x})-L_x(t,x,\dot{x})=0,    \label{2.7}\\
&& x(0)=x(T),  \qquad \dot{x}(0)=\dot{x}(T).    \label{2.8}\eea

For such an extremal loop, define
\bea
P(t) &=& L_{p,p}(t,x(t),\dot{x}(t)),  \nn\\
Q(t) &=& L_{x,p}(t,x(t),\dot{x}(t)),  \nn\\
R(t) &=& L_{x,x}(t,x(t),\dot{x}(t)).  \nn\eea
Note that
\be F\,''(x)=-\frac{d}{dt}(P\frac{d}{dt}+Q)+Q^T\frac{d}{dt}+R. \lb{2.9}\ee

For $\omega\in\U$, set
\be  D(\omega,T)=\{y\in W^{1,2}([0,T],\C^n)\,|\, y(T)=\omega y(0) \}.   \lb{2.10}\ee
We define the $\omega$-Morse index $\phi_\omega(x)$ of $x$ to be the dimension of the
largest negative definite subspace of
$$ \langle F\,''(x)y_1,y_2 \rangle, \qquad \forall\;y_1,y_2\in D(\omega,T), $$
where $\langle\cdot,\cdot\rangle$ is the inner product in $L^2$. For $\omega\in\U$, we
also set
\be  \ol{D}(\omega,T)= \{y\in W^{2,2}([0,T],\C^n)\,|\, y(T)=\omega y(0), \dot{y}(T)=\om\dot{y}(0) \}.
                     \lb{2.11}\ee
Then $F''(x)$ is a self-adjoint operator on $L^2([0,T],\R^n)$ with domain $\ol{D}(\omega,T)$.
We also define
\[\nu_\omega(x)=\dim\ker(F''(x)).\]

In general, for a self-adjoint operator $A$ on the Hilbert space $\mathscr{H}$, we set
$\nu(A)=\dim\ker(A)$ and denote by $\phi(A)$ its Morse index which is the maximum dimension
of the negative definite subspace of the symmetric form $\langle A\cdot,\cdot\rangle$. Note
that the Morse index of $A$  is equal to the total multiplicity of the negative eigenvalues
of $A$.

On the other hand, $\td{x}(t)=(\partial L/\partial\dot{x}(t),x(t))^T$ is the solution of the
corresponding Hamiltonian system of (\ref{2.7})-(\ref{2.8}), and its fundamental solution
$\gamma(t)$ is given by
\bea \dot{\gamma}(t) &=& JB(t)\gamma(t),  \lb{2.12}\\
     \gamma(0) &=& I_{2n},  \lb{2.13}\eea
with
\be B(t)=\left(\matrix{P^{-1}(t)& -P^{-1}(t)Q(t)\cr
                       -Q(t)^TP^{-1}(t)& Q(t)^TP^{-1}(t)Q(t)-R(t)\cr}\right). \lb{2.14}\ee

\begin{lemma}(Long, \cite{Lon4}, p.172)\lb{L2.3}  
For the $\omega$-Morse index $\phi_\omega(x)$ and nullity $\nu_\omega(x)$ of the solution $x=x(t)$
and the $\omega$-Maslov-type index $i_\omega(\gamma)$ and nullity $\nu_\omega(\gamma)$ of the
symplectic path $\ga$ corresponding to $\td{x}$, for any $\omega\in\U$ we have
\be \phi_\omega(x) = i_\omega(\gamma), \qquad \nu_\omega(x) = \nu_\omega(\gamma).  \lb{2.15}\ee
\end{lemma}

A generalization of the above lemma to arbitrary  boundary conditions is given in \cite{HS1}.
For more information on these topics, we refer to \cite{Lon4}.

\medskip

\noindent {\bf Acknowledgements.} The authors thank sincerely Professor Xijun Hu, especially for
discussions with him on the proof of Theorem \ref{Th:no.intersection}. They thank the anonymous
referee sincerely for his/her many helpful suggestions and comments on the first manuscript of
this paper.


\begin{thebibliography}{}

\bibitem{Dan} J. Danby, The stability of the triangular
Lagrangian point in the general problem of three bodies.
{\it Astron. J.} 69. (1964) 294-296.

\bibitem{Euler} L. Euler, De motu restilineo trium corporum se mutus
attrahentium. {\it Novi Comm. Acad. Sci. Imp. Petrop.}  11. (1767) 144-151.

\bibitem{Ga} M. Gascheau,  Examen d'une classe d'\'{e}quations diff\'{e}rentielles
et application \`{a} un cas particulier du probl\`{e}me des trois corps. {\it Comptes
Rend. Acad. Sciences.} 16. (1843) 393-394.

\bibitem{Go} W. B. Gordon, A minimizing property of Kepler
orbits, {\it American J. of Math.} 99. (1977) 961-971.

\bibitem{HLS} X. Hu, Y. Long, S. Sun, Linear stability of elliptic Lagrange
solutions of the classical planar three-body problem via index theory.
{\it Arch. Ration. Mech. Anal.} 213. (2014) 993-1045.

\bibitem{HO} X. Hu, Y. Ou,  An estimate for the hyperbolic region of elliptic
Lagrangian solutions in the planar three-body problem.
{\it Regul. Chaotic. Dyn.} 18(6). (2013) 732-741.

\bibitem{HO2} X. Hu, Y. Ou,  Collision index and stability of elliptic relative
equilibria in planar $n$-body problem. {\it https://arxiv.org/pdf/1509.02605.pdf},
(2015). {\it Comm. Math. Phys.} to appear.

\bibitem{HS1} X. Hu, S. Sun,  Index and stability of symmetric periodic
orbits in Hamiltonian systems with its application to figure-eight orbit.
{\it Commun. Math. Phys.} 290. (2009) 737-777.

\bibitem{HS2} X. Hu, S. Sun, Morse index and stability of elliptic Lagrangian
solutions in the planar three-body problem. {\it Advances in Math.} 223. (2010) 98-119.

\bibitem{Jac1} N. Jacobson, Basic Algebra I. W. H. Freeman and Com. 1974.

\bibitem{Ka} T. Kato, Perturbation Theory for Linear Operators. Second
edition, Springer-Verlag, Berlin, 1984.

\bibitem{Lag} J. Lagrange,  Essai sur le probl\`{e}me des
trois corps. Chapitre II. {\OE}uvres Tome {6}, Gauthier-Villars,
Paris. (1772) 272-292.

\bibitem{Lon1} Y. Long, Maslov-type index, degenerate critical points, and
asymptotically linear Hamiltonian systems. {\it Science in China}.
Series A. 7. (1990) 673-682. (Chinese Ed.). Series A. 33. (1990)
1409-1419. (English Ed.).

\bibitem{Lon2} Y. Long, Bott formula of the Maslov-type index
theory. {\it Pacific J. Math.} 187. (1999) 113-149.

\bibitem{Lon3} Y. Long, Precise iteration formulae of the Maslov-type
index theory and ellipticity of closed characteristics. {\it Advances
in Math.} 154. (2000) 76-131.

\bibitem{Lon4} Y. Long, Index Theory for Symplectic Paths with
Applications. Progress in Math. 207, Birkh\"auser. Basel. 2002.

\bibitem{Lon5} Y. Long, Lectures on Celestial Mechanics and Variational Methods.
{\it Preprint.} 2012

\bibitem{MaS1} R. Mart\'{i}nez, A. Sam$\grave{a}$, On the centre mabifold of
collinear points in the planar three-body problem. {it Cele. Mech. and Dyn. Astro.}
85. (2003) 311-340.

\bibitem{MSS} R. Mart\'{\i}nez, A. Sam\`{a}, C. Sim\'{o},
Stability of homograpgic solutions of the planar three-body problem
with homogeneous potentials. in International conference on
Differential equations. Hasselt, 2003, eds, Dumortier, Broer, Mawhin,
Vanderbauwhede and Lunel, World Scientific, (2004) 1005-1010.

\bibitem{MSS1} R. Mart\'{\i}nez, A. Sam\`{a}, C. Sim\'{o},
Stability diagram for 4D linear periodic systems with applications
to homographic solutions. {\it J. Diff. Equa.} 226. (2006) 619-651.

\bibitem{MSS2} R. Mart\'{\i}nez, A. Sam\`{a}, C. Sim\'{o}, Analysis of
the stability of a family of singular-limit linear periodic systems in
$\R^4$. Applications. {\it J. Diff. Equa.} 226. (2006) 652-686.

\bibitem{MS} K. Meyer, D. Schmidt, Elliptic relative equilibria in
the N-body problem. {\it J. Diff. Equa.} 214. (2005) 256-298.

\bibitem{R1} G. Roberts, Linear stability of the elliptic Lagrangian
triangle solutions in the three-body problem. {\it J. Diff. Equa.}
182. (2002) 191-218.

\bibitem{R2} E. Routh, On Laplace's three particles with
a supplement on the stability or their motion. {\it Proc. London Math.
Soc.} 6. (1875) 86-97.

\bibitem{V1} A. Venturelli, Une caract\'erisation variationelle
des solutions de Lagrange du probl\'eme plan des trois corps. {\it C. R.
Acad. Sci. Paris S\'er.} I. 332. (2001) 641-644.

\bibitem{Win1} A. Wintner, The Analytical Foundations of Celestial Mechanics.
Princeton Univ. Press, Princeton, NJ. 1941. Second print, Princeton Math. Series
5, 215. 1947.

\bibitem{ZZ1} S. Zhang, Q. Zhou, A minimizing property of Lagrangian
solutions. {\it Acta Math. Sin. (Engl. Ser.)} 17. (2001) 497-500.

\bibitem{ZhL1} G. Zhu, Y. Long, Linear stability of some symplectic matrices.
{\it Frontiers of Math. in China.} 5. (2010) 361-368.

\bibitem{ZhoL1} Q. Zhou, Y. Long, Maslov-type indices and linear stability
of elliptic Euler solutions of the three-body problem. {\it http://arxiv.org/abs/1510.06822.pdf.
} (2015).

\bibitem{ZhoL2} Q. Zhou, Y. Long, The reduction on the linear stability of elliptic
Euler-Moulton solutions of the $n$-body problem to those of $3$-body problems.
{\it http://arxiv.org/abs/1511.00070.pdf.} (2016), submitted.

\end{thebibliography}
\end{document}